\documentclass[11pt]{article}

\usepackage{epsf}
\usepackage{psfrag}
\usepackage{amsmath}
\usepackage{amssymb}
\usepackage{lscape}
\usepackage{latexsym}
\usepackage[usenames]{color}
\usepackage[mathcal]{eucal}
\usepackage{stmaryrd}
\usepackage{textcomp}
\usepackage{setspace} 
\usepackage{graphicx} 
\usepackage{pgfplots} 
\usepackage{epsfig}   
\usepackage{epstopdf} 
\usepackage{appendix}
\usepackage{afterpage}
\usepackage{float} 
\usepackage{array} 
\usepackage[misc,geometry]{ifsym} 
\usepackage[font=small,skip=0pt]{caption}

\newtheorem{theorem}{Theorem}[section]

\newtheorem{remark}{Remark}[section]

\newtheorem{algorithm}{Algorithm}[section]

\newcommand{\ds}{\displaystyle}
\newcommand{\dZ}{{\cal Z \kern -0.7em Z}}
\newcommand{\dC}{{\rm\hbox{C \kern-0.8em\raise0.2ex\hbox{\vrule height5.4pt width0.7pt}}}}
\newcommand{\dQ}{{\rm\hbox{Q \kern-0.85em\raise0.25ex\hbox{\vrule height5.4pt width0.7pt}}}}

\newcommand\old[1]{}
\newcommand{\beqa}{\begin{eqnarray*}}
\newcommand{\eeqa}{\end{eqnarray*}}

\title{\bf \Large{{A Path Planning Model for Intercepting a Moving Target with Finite Obstacle Avoidance}}}

\author{Masuda Akter  \thanks{Department of Mathematics, University of Chittagong, Chittagong-4331, Bangladesh,
   {\tt email: masuda@cu.ac.bd}.} 
\and
M. M. Rizvi \Letter \thanks{Centre for Smart Analytics (CSA), Institute of Innovation, Science and Sustainability, Federation University Australia, {\tt email:  rizmm001@mymail.unisa.edu.au}, (corresponding author).}}
\begin{document}

\maketitle

\begin{abstract}
  \noindent This paper investigates the problem of computing a two-dimensional curvature-constrained shortest path for an unmanned aerial vehicle (UAV) to intercept a moving target, with both the UAV (pursuer) and the target traveling at constant speeds. We formulate an optimal control problem that integrates two critical objectives: avoiding static obstacles and successfully intercepting the target. The formulation integrates curvature constraints, geometric obstacle-avoidance conditions, and time-consistent interception constraints within a constrained optimal control framework. A geometric framework is developed, along with necessary conditions for path optimality under the imposed constraints. The problem is initially examined in the presence of a single obstacle and later extended to scenarios involving a finite number of obstacles. Numerical experiments are carried out to evaluate the performance and efficiency of the proposed model using illustrative examples. Finally, we present a realistic case study using actual geographic data, including obstacle placement, target trajectory, and heading angles, to demonstrate the practical applicability and effectiveness of the proposed method in real-world scenarios.

\end{abstract}


\textbf{AMS} subject classifications: 34K35 49M37, 90C30, 	90C39.

\section{Introduction}\label{sec1}
\noindent The development of optimal path planning techniques is crucial to improving movement efficiency, reducing operational costs, and ensuring safety and convenience in various real-world applications. Path planning is especially essential in fields such as robotics \cite{Yong2006}, agriculture \cite{Höffmann2024}, and rescue operations \cite{chou2019}. A critical challenge in autonomous systems 
and robotic navigation involves devising strategies to intercept a moving target while avoiding static obstacles positioned between the pursuer and the target. The primary objective of this study is to formulate an optimal path-planning model that takes into account the presence of static obstacles, allowing the pursuer to navigate around them while intercepting the moving target efficiently. This research seeks to contribute to the theoretical foundations and practical applications of path planning strategies in dynamic environments by employing mathematical modelling and efficient algorithms. Through the development of this model, we aim to address the complexities of intercepting a moving target amid various obstacles, providing a model applicable across a range of domains.

\noindent The shortest path was first proposed by Dubins \cite{Dubins1957} for vehicles that only travel forward and make left or right unidirectional turns from predetermined starting and finishing places in the 2D plane. The author defines paths as combinations of straight ($S$) and curves ($C$) lines, where the vehicle may be straight or curved in the middle and starts and finishes with turns. $Z=\{LSL, RSL, RSR, LSR, LRL, RLR\}$ represents the six potential permutations of these pathways, where L and R stand for left and right turns, respectively. Later on, Dubins work was extended by Reeds and Shepp \cite{Reeds1990}, allowing for backward adaptability in these routes. Since then, Dubins findings have been expanded upon, utilizing sophisticated theories in nonlinear optimal control \cite{Sussmann1991}-\cite{Chen2019} and integrating geometric concepts \cite{Ayala2014}. In particular, the application of Pontryagin's Maximum Principle (PMP) \cite{Pontryagin1962} has improved knowledge and uses of Dubins pathways. Due to these advancements, Dubins work now encompasses a broader range and is considered a fundamental basis for the research of optimized vehicle trajectories. 

\noindent Unmanned aerial vehicles (UAVs) often carry out surveillance missions when they are directed to visit specific target sites and gather the necessary data. The Dubins Touring Problem (DTP) \cite{Cohen2017} can be used to frame the problem if the order of visits to the destinations is known in advance.   Many studies \cite{Madveev2011, Madveev2012} have been conducted on Dubins path problem for intercepting moving targets, looking at both relaxed optimal pathways and Dubins \cite{Gopalan2016}-\cite{Buzikov2021}. Looker \cite{Looker2008} suggests a search algorithm to find the shortest route. Meyer et al. \cite{ Meyer2015} study minimum-time interception, predetermined time for target movement, and the relaxed Dubins problem, focusing on finding the shortest path for a Dubins vehicle without a terminal angle constraint. In \cite{Yang2002}, Yang and Kapila investigate two-dimensional optimal path-planning problems for UAVs, accounting for both kinematic and tactical constraints. By applying vector calculus, the problem is reformulated as a parameter optimization task. An efficient numerical algorithm is developed to satisfy the optimality conditions. The proposed approach can manage multiple tactical constraints simultaneously, and its effectiveness is demonstrated through a representative numerical simulation.
The optimality principles for collision-free trajectory and Dubins interpolating curves—the shortest curvature-constrained pathways in a given sequence—are the main topics of Kaya's research \cite{Kaya2017,Kaya2019} on Dubins path reformulation. Zheng et al. introduced a model for the Minimum-Time Intercept Problem \cite{Zheng2021}, guiding a Dubins vehicle to intercept a moving target with prescribed impact angle constraints \cite{Zheng2022}. Forkan et al. \cite{Forkan2022}  developed a mathematical model using arc parameterization to determine the optimal path length or time for touring a finite number of targets using an unmanned aerial vehicle. A method was designed for deploying a mobile sensor network to detect and capture moving targets in a planar environment in \cite{ Ferrari2009}. In order to guarantee target capturability and zero terminal guidance command, Hou \cite{Hou2023} suggests geometry-based impact time and angle control rules that do not require time-to-go calculations.

\noindent The significance of UAV path planning for mission performance has led to the development of numerous approaches, particularly for challenges involving obstacle avoidance restrictions  \cite{Cao2015}. 
Recently, the shortest Dubins path between three consecutive via-points with given initial and final orientations is found using a novel technique presented in the research \cite{Parlangeli2024} that makes use of analytic geometry.

\noindent In this research, we formulate an optimal control problem to determine the shortest interception path for a pursuer aiming to reach a moving target while avoiding obstacles. We develop mathematical expressions to calculate the lengths of individual circular and straight line sub paths and construct a model to account for the presence of obstacles. In our model, the position of the pursuer within the obstacle’s turning circle is dynamically determined as a function of time. Unlike previous studies that considered fixed target points in path planning, our research introduces a framework in which both the touring point and the interception point are unknown. This makes a significant contribution to the field of path planning by enabling more realistic and adaptive interception scenarios. Several solvers are employed to solve the problem, as discussed in the numerical analysis section. Additionally, a real-world case study is presented to test the model's capability and evaluate its effectiveness in computing an optimal interception path.

\noindent In Section \ref{SM}, we define our problem statement and the motivation behind the study. The main contribution and outcome are briefly described in Section \ref{CO}. Section \ref{geometry1} presents the geometric structure of the engagement of the pursuer, obstacle, and target. Building on this geometry, Section \ref{MT} formulates the optimal control problem and characterizes the optimal control of the pathway. The proposed path planning model for intercepting a moving target touring with static obstacles is demonstrated in Section \ref{modelM}. Section \ref{NM} discusses the findings and interpretation of numerical tests for different scenarios. Later in Section \ref{RWA}, we present a real-life scenario of our proposed model. The conclusion of the paper is presented in the concluding part.

\section{ Problem Statement and Motivation}\label{SM}
\noindent This study addresses a fundamental problem in real-world scenarios where a pursuer, moving at a constant speed, aims to intercept a moving target, also traveling at a constant speed, while navigating around a static obstacle situated between the pursuer and the target. The target is assumed to move in a straight line. This problem is particularly challenging in the contexts of autonomous systems, robotics, and unmanned vehicles, where the goal is to intercept a moving target while avoiding obstacles that may obstruct the pursuer's path. Real-world applications, such as autonomous drone navigation, search and rescue missions, and surveillance operations, often present similar challenges. In these scenarios, the pursuer must reach the target quickly while effectively maneuvering around obstacles that may be positioned between them. While most conventional path planning models prioritize direct interception or obstacle avoidance, there has been limited research on the concept of "touring" that is, circumventing obstacles in a more strategic manner to approach the target. This gap motivates the present study, which seeks to develop a unified framework capable of simultaneously handling obstacle avoidance and moving target interception under minimum-path-length criteria.

\noindent The motivation for this work stems from the need to develop more realistic and computationally tractable models for UAV path planning. Existing approaches either simplify the obstacle environment, assume static target positions, or do not adequately capture the curvature constraints inherent to fixed-wing UAV dynamics. By integrating Dubins path theory with a dynamic interception model, this paper aims to fill these gaps and provide a model that is both theoretically sound and practically applicable.

\noindent In the numerical experiments Section \ref{NM}, we demonstrate that the proposed approach achieves a more flexible and realistic solution for the interception scenario. These results indicate that the approach can support interception planning in environments with static obstacles.
\section{Contribution and Outcome}\label{CO}



This section summarizes the principal theoretical and computational contributions of the present work:
\begin{enumerate}
    \item [(i)] Optimal Control Formulation: We formulate a constrained optimal control problem for determining the minimum-length curvature-constrained path (CS path) that enables a pursuer (UAV) to intercept a moving target while avoiding one or more static circular obstacles.
 \item [(ii)] Dynamic Touring Point: Unlike prior studies that assume fixed waypoints, the touring point, the point at which the pursuer enters the turning circle of the obstacle, is treated as an unknown variable determined dynamically as a function of time and the relative geometry of the pursuer and target.

 \item [(iii)] Unknown Interception Point: The interception point is not prescribed in advance. Instead, it is computed as part of the optimization, allowing the framework to adapt to varying target speeds, headings, and initial conditions.
 
 \item [(iv)] Geometric Optimality Conditions: A set of necessary geometric conditions is derived to characterize optimal CS paths under the imposed obstacle avoidance and interception constraints, enabling efficient computation by restricting the search to four candidate path types.

 \item [(v)] Extension to Multiple Obstacles: The model is systematically extended from a single-obstacle scenario to configurations involving a finite number of static obstacles, with a general algorithmic strategy for sequentially resolving each obstacle along the interception path.
  
   \item [(vii)] Real-World Case Study: A realistic scenario using actual geographic data  including terrain-based obstacle placement, target trajectory, and heading angles is presented to validate the proposed model's practical applicability.
\end{enumerate}
\section {Geometric Pathway of Interception} \label{geometry1}


\begin{figure}[!htbp]
\hspace{-1cm}
\begin{minipage}{80mm}
\begin{center}
\hspace*{0cm}
\includegraphics[width=70mm]{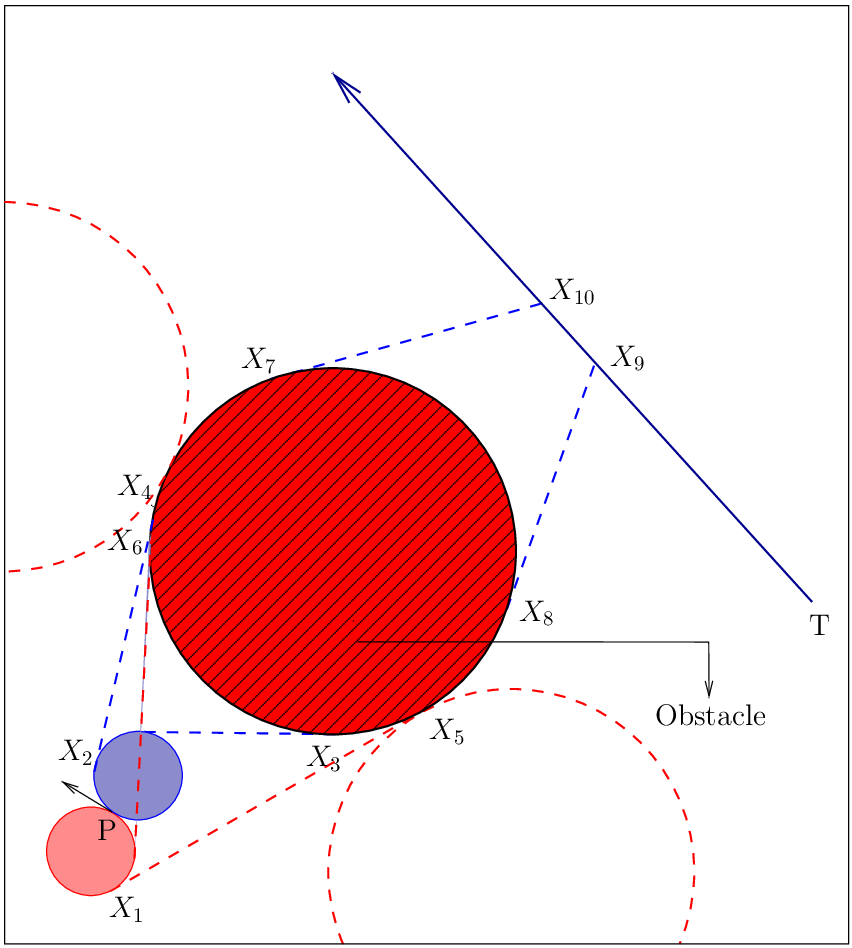} \\ \hspace{-3mm}
{\footnotesize (a) Circular Obstacle.}\hspace{-3mm}
\end{center}
\end{minipage}
\hspace{-6mm}
\begin{minipage}{80mm}
\begin{center}
\hspace*{0cm}
\includegraphics[width=70mm]{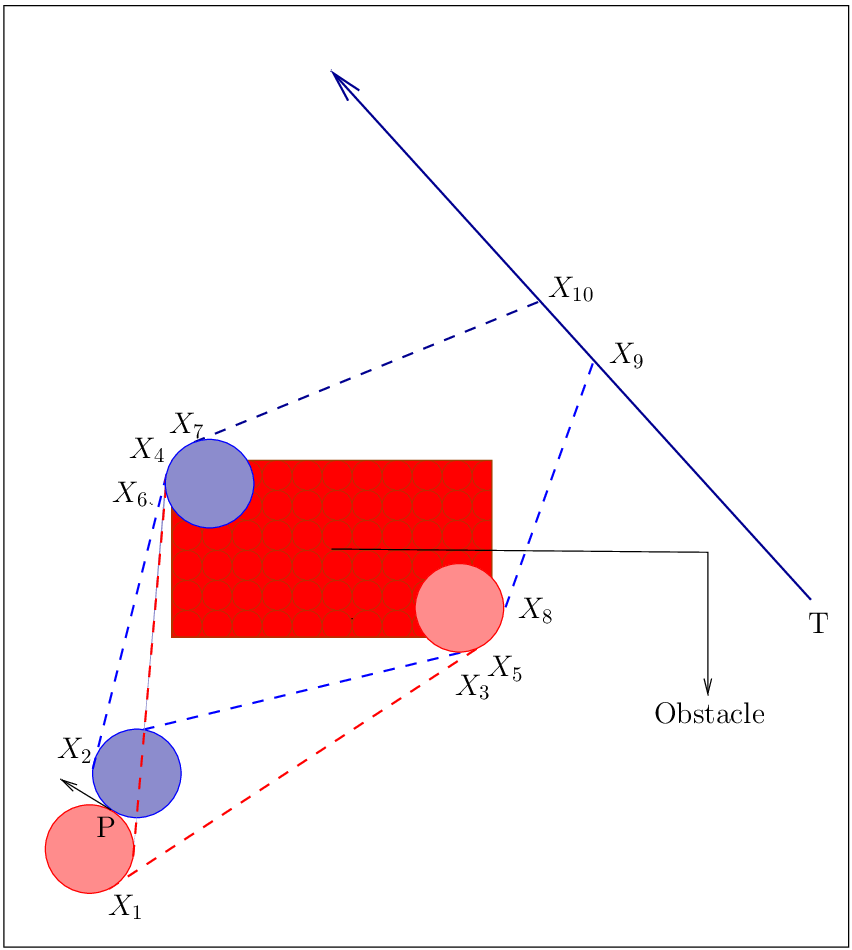} \\ \hspace{-2mm}
{\footnotesize (b) Rectangular Obstacle.}\hspace{-2mm}
\end{center}
\end{minipage}
\\\hspace{-2mm}

\vspace{.5cm}
\caption{\sf\small Paths of interception of a moving target and a pursuer touring a single obstacle. }
\label{Figopti_path}
\end{figure}

\noindent This section presents the geometric representation of the interception path between a pursuer and a moving target. The geometric decomposition introduced here is used in Section \ref{MT} to state the optimal control problem: every circular arc of the pursuer's turning circle corresponds to a bang control $u(t)=\pm a$, every straight segment to $u(t)=0$, and every arc along the obstacle turning circle to a boundary arc of the obstacle-avoidance constraint, on which the curvature is $w(t)=\pm a_w$.  It is assumed that the initial positions and heading angles of both the pursuer $P(x_P(t_0), y_P(t_0), \theta_P(t_0))$ and the target $T(x_T(t_0), y_T(t_0), \theta_T(t_0))$ are known. In environments with static obstacles, the center and radius of each obstacle are also assumed to be known in advance. When an obstacle has a rectangular shape, a circular arc is introduced at the point where the pursuer begins a turning maneuver. The radius of this arc is set equal to the minimum turning radius of the pursuer to ensure the path remains feasible. This geometric approach provides both a visual and analytical understanding of the interception process in environments with obstacles.
\color{black}

\noindent According to Figure \ref{Figopti_path}, the pursuer follows either path $PX_1$ or $PX_2$ to travel the circular path of the pursuer turning circle. The pursuer then follows one of the four straightline $X_1X_5$, $X_1X_6$, $X_2X_3$, $X_2X_4$ to enter the obstacle turning circle. Then the pursuer tours the obstacle in one of the four feasible ways $X_3X_8$, $X_5X_8$, $X_4X_7$, $X_6X_7$. After touring the obstacle, the pursuer follows the straight way to intercept the target, either $X_7X_{10}$ or $X_8X_9$. This creates four possible feasible paths of the type $CS+CS+S_T$ which are as follows:
 \begin{equation}\label{geo_inter}
  \footnotesize{ CS+CS+S_T\equiv \\     
   \begin{cases}
   LS +LS + S_T\equiv \\ \mbox{arc}\;PX_1 + \mbox{st. line}\;X_1X_5 + \mbox{arc}\;X_5X_8+\mbox{st. line}\;X_8X_9 + \mbox{st. line}\;TX_9,&\text {or},\\
   LS +RS + S_T\equiv\\ \mbox{arc}\;PX_1 + \mbox{st. line}\;X_1X_6 + \mbox{arc}\;X_6X_7+\mbox{st. line}\;X_7X_{10} + \mbox{st. line}\;TX_{10},&\text {or},\\
   RS +LS + S_T\equiv\\ \mbox{arc}\;PX_2 + \mbox{st. line}\;X_2X_3 + \mbox{arc}\;X_3X_8+\mbox{st. line}\;X_8X_9 + \mbox{st. line}\;TX_9,&\text {or},\\
   RS +RS + S_T\equiv \\ \mbox{arc}\;PX_2 + \mbox{st. line}\;X_2X_4 + \mbox{arc}\;X_4X_7+\mbox{st. line}\;X_7X_{10} + \mbox{st. line}\;TX_{10},&\text {or},\\
  \end{cases} }
\end{equation}
\noindent Here, $L$ and $R$ represent the left and right circular arcs, respectively, while $S$ denotes the straight-line segment of the pursuer's path. Additionally, $S_T$ represents the straight-line path of the target.

\section{Optimal Control Formulation}\label{MT}
\noindent In this section, we formulate the optimal control problem associated with the geometric pathway described in Section \ref{geometry1} and characterize the structure of its optimal control. The pursuer can successfully intercept the target only if its speed exceeds that of the target. We consider the problem in two-dimensional space.

\noindent 
Consider a pursuer moving in a two-dimensional plane, with its position at time \( t \) given by \((x_p(t), y_p(t)) \in \mathbb{R}^2\). The pursuer follows an optimal trajectory \(\ell(t): [t_0, t_f] \to \mathbb{R}^2\), aiming to intercept a moving target while avoiding fixed obstacles along the way. We assume that the target moves along a straight-line path with a constant velocity, and the pursuer has complete knowledge of the target’s motion, including its speed and direction.  

\noindent A stationary obstacle, centered at \((x_b, y_b)\), is positioned along the potential interception path, requiring the pursuer to maneuver around it while maintaining an effective pursuit strategy. The pursuer's speed is denoted by \(|\mathbf{V_p}|\), while the target's speed is \(|\mathbf{V_T}|\). Additionally, the pursuer is subject to a minimum turning radius constraint, denoted by \( R \), which limits its ability to make sharp turns; equivalently, the curvature of its path is bounded by $a=1/R$. The obstacle is enclosed by a circle, called the obstacle turning circle, of radius $R_b=1/a_w$ centered at $(x_b,y_b)$, which the pursuer must not enter. We assume $R_b\geq R$ (that is, $a_w\leq a$), so that the pursuer is able to fly along the obstacle turning circle; an obstacle with sharp corners is enclosed by such a circle as illustrated in Figure \ref{Figopti_path}(b). This constraint necessitates an optimal trajectory that balances efficient interception with obstacle avoidance, ensuring the pursuer reaches the target in the shortest possible time. Let the kinematic equations of the pursuer's path be
\begin{equation}\label{kinematic}
\begin{bmatrix}
\dot{x_p}(t)  \\
\dot{y_p}(t) \\
\dot{\theta_p}(t)
\end{bmatrix}
=
\begin{bmatrix}
|\mathbf{V_p}|\cos\theta_P(t)  \\
|\mathbf{V_p}|\sin\theta_p(t) \\
|\mathbf{V_p}|u(t)
\end{bmatrix},   
\end{equation} 
\noindent where $\theta_p(t)\in [0,2\pi)$ denotes the heading angle of the pursuer, measured counterclockwise from the positive $x$-axis.
Here the curvature is the main control variable $u(t)$. The control is calculated as $u(t) = \ds \frac{\dot{\theta_p}(t)}{|\mathbf{V_p}|}$, which is positive or negative depending on whether the pursuer takes a left or right turn, respectively, and zero when the pursuer moves in a straight line. The total arc length travelled by the pursuer,
\begin{equation}\label{minilengthpath}
\int_{t_0}^{t_f} \sqrt{\dot{x_p}^2+\dot{y_p}^2}\,dt = \int_{t_0}^{t_f} ||\dot{\ell}(t)||\,dt=|\mathbf{V_p}|(t_f-t_0),   
\end{equation}
where $||\dot{\ell}(t)||=|\mathbf{V_p}|$, $\ds \dot{\ell}(t)=\frac{d\ell}{dt}$, and $||\cdot||$ is the Euclidean norm, is adopted as the objective (cost) functional of the optimal control problem \eqref{problem1} below. The functional \eqref{minilengthpath} by itself only measures the length of a given path; a minimal-length path is obtained only when \eqref{minilengthpath} is minimized subject to the pursuer's kinematics \eqref{kinematic}, the curvature bound, the obstacle-avoidance constraint and the interception conditions stated in \eqref{problem1}. The resulting trajectory is the pursuer's minimal-length path for touring a static obstacle and intercepting a moving target.

 \noindent When the pursuer navigates around a static obstacle, it either follows its minimum turning radius if the obstacle has a sharp corner or adopts the obstacle's turning radius. Consequently, the curvature of the pursuer's path changes, and this new curvature is denoted by \( w(t) \).
 It is assumed that the kinematic equations \eqref{pursuerkinematic} describing the pursuer's motion while touring around the obstacle are given as follows.
\begin{equation}\label{pursuerkinematic}
\begin{bmatrix}
\dot{x_p}(t)  \\
\dot{y_p}(t) \\
\dot{\theta_b}(t)
\end{bmatrix}
=
\begin{bmatrix}
|\mathbf{V_p}|\cos\theta_b(t)  \\
|\mathbf{V_p}|\sin\theta_b(t) \\
|\mathbf{V_p}|w(t)
\end{bmatrix},   
\end{equation}  
where $w(t) = \ds \frac{\dot{\theta_b}(t)}{|\mathbf{V_p}|}$, and $\theta_b(t)$ denotes the heading angle of the pursuer while it moves along the obstacle turning circle. Thus, $w(t)$ is not an additional control associated with the stationary obstacle; it is the curvature $u(t)$ of the pursuer itself on those parts of the path where the pursuer flies along the boundary of the obstacle turning circle.

\noindent Suppose that the target location is ($x_T(t)$, $y_T(t)$, $\theta_T$) and that it moves without turns or maneuvers at a constant speed $|\mathbf{V_T}|>0$. Within $[0, 2\pi]$, $\theta_T$ indicates the direction of the target's movement. Throughout the engagement, the value of $\theta_T$ stays constant because the target does not change its orientation. Considering that the target's starting position is $(x_{T_0},y_{T_0})$, the following equations represent its motion:

\begin{equation}\label{problemT}
\begin{array}{cl} 
   &\  \dot{x}_T(t)=|\mathbf{V_T}| \cos\theta_T, \,x_T(t_0)= x_{T_0}, \\[2mm]
    &  \ \dot{y}_T(t)= |\mathbf{V_T}|\sin\theta_T,\, y_T(t_0)= y_{T_0}, \\[2mm]
\end{array}\\
\end{equation}
whose solution is

\begin{equation}\label{targetsol}
x_T(t)=x_{T_0}+|\mathbf{V_T}|\cos\theta_T\,(t-t_0), \qquad y_T(t)=y_{T_0}+|\mathbf{V_T}|\sin\theta_T\,(t-t_0).
\end{equation}

\noindent Therefore, the optimal control problem of intercepting the moving target \eqref{problemT} while avoiding the obstacle can be formulated as follows.

\begin{equation}\label{problem1}
\begin{array}{cll} \ds\min_{u(\cdot),\;t_f} & \ f =\ds \int_{t_0}^{t_f} |\mathbf{V_P}| \,dt =|\mathbf{V_P}| (t_f-t_0) & \\[4mm]
\mbox{subject to} & & \\[1mm]
&\  \dot{x}_P(t)=|\mathbf{V_P}| \cos\theta_P(t), & x_P(t_0)= x_{P_0}, \\[2mm]
&  \ \dot{y}_P(t)= |\mathbf{V_P}|\sin{\theta}_P(t), & y_P(t_0)= y_{P_0} , \\[2mm]
& \ \dot{\theta}_P(t) =u(t)|\mathbf{V_P}|, & \theta_P(t_0)=\theta_{P_0},  \\[2mm]
&  \ |u(t)|\leq a, & t\in[t_0,t_f],\\[2mm]
&  \ \big(x_P(t)-x_b\big)^2+\big(y_P(t)-y_b\big)^2\geq R_b^2, & t\in[t_0,t_f],\\[2mm]
&  \ x_P(t_f)=x_T(t_f),\quad y_P(t_f)=y_T(t_f). &
\end{array}
\end{equation}

\noindent where ($x_P(t)$, $y_P(t)$, $\theta_P(t)$) represents the position and heading along the path of the pursuer, and $(x_T(t_f),y_T(t_f))$ is given by \eqref{targetsol}. The decision variables of \eqref{problem1} are the curvature control $u(\cdot)$ and the free terminal (interception) time $t_f$; the state trajectory $(x_P,y_P,\theta_P)$ is determined by $u(\cdot)$ through the kinematics and the initial conditions. The terminal conditions $x_P(t_f)=x_T(t_f)$ and $y_P(t_f)=y_T(t_f)$ require the pursuer and the target to occupy the same position at the interception time and thereby link the target's motion to the pursuer's objective. No terminal condition is imposed on the heading $\theta_P(t_f)$, since only position interception, without a prescribed impact angle, is required; the pursuer may therefore arrive at the interception point with any heading, which is determined by the optimization together with the interception time $t_f$. The pure state constraint $(x_P(t)-x_b)^2+(y_P(t)-y_b)^2\geq R_b^2$ keeps the pursuer outside the obstacle turning circle for all $t\in[t_0,t_f]$. For $n$ obstacles with centres $(x_b^j,y_b^j)$ and radii $R_b^j$, this constraint is imposed for each $j=1,\ldots,n$.

\noindent On a subarc where the obstacle-avoidance constraint is active, the pursuer moves along the obstacle turning circle. Differentiating the active constraint along the trajectory shows that the curvature on such a boundary arc is constant and equal to $w(t)=+a_w$ (counterclockwise touring, obstacle on the pursuer's left) or $w(t)=-a_w$ (clockwise touring, obstacle on the pursuer's right), with $|w(t)|\leq a_w\leq a$, so the curvature bound $|u(t)|\le a$ is respected. Hence, the obstacle turning circle with control $w(t)$ introduced in \eqref{pursuerkinematic} describes exactly the boundary arcs of the state constraint in \eqref{problem1}, and the pursuer enters and leaves such an arc tangentially. If the obstacle does not obstruct the interception path, the obstacle-avoidance constraint is inactive and \eqref{problem1} reduces to the obstacle-free interception problem.

\noindent To generate the shortest path for the interception of a moving target, the model is expected to first identify all feasible trajectories that the pursuer can follow. These trajectories are determined based on the pursuer’s ability to turn left or right, considering both its minimum turning radius and the obstacle’s turning radius.  After generating this set of candidate paths, the model selects the optimal trajectory that minimizes both the total path length and the time required to intercept the moving target.

\subsection{Optimality Conditions via Pontryagin’s Maximum Principle}

In this section, we derive necessary conditions for optimality for the interception problem formulated in \eqref{problem1} and \eqref{problemT} using Pontryagin’s Maximum Principle (PMP) (see \cite{Pontryagin1962}). The objective is to characterize the structure of the optimal control and establish the geometric form of optimal trajectories. The analysis below applies to the interior subarcs of an optimal trajectory, that is, to the parts on which the obstacle-avoidance constraint in \eqref{problem1} is inactive; on boundary subarcs the curvature is $w(t)=\pm a_w$, as discussed above.

\noindent Since the pursuer’s speed (\(|V_P|\)) is constant, minimizing the function defined in \eqref{problem1} is equivalent to minimizing the terminal time (\(t_f\)). Therefore, the problem can be reformulated as a minimum-time optimal control problem with bounded curvature control.

\noindent Hamiltonian Formulation: Let the state vector of the pursuer be defined as
\begin{equation}\label{statevector}
x(t) = (x_P(t), y_P(t), \theta_P(t))^\top,
\end{equation}
with dynamics given in \eqref{problem1}, and control ($u(t) \in [-a,a]$). The Hamiltonian associated with this problem is defined in \eqref{hamiltonian}.
\begin{equation}\label{hamiltonian}
    H(x,\lambda,u) = 1 + \lambda_1 |V_P| \cos\theta_P + \lambda_2 |V_P| \sin\theta_P
+ \lambda_3 |V_P| u,
\end{equation}

\noindent where ($\lambda(t) = (\lambda_1, \lambda_2, \lambda_3)^\top$) denotes the adjoint (costate) vector.

\noindent According to PMP, the costate variables satisfy
\begin{equation}\label{costate}
    \dot{\lambda}(t) = -\frac{\partial H}{\partial x}(x^*(t),\lambda(t),u^*(t)).
\end{equation}

\noindent Computing partial derivatives of \eqref{hamiltonian} yields the adjoint system
\begin{equation}\label{adjoint}
\dot{\lambda}_1 = 0, \qquad \dot{\lambda}_2 = 0, \qquad \dot{\lambda}_3 = -|V_P|(-\lambda_1 \sin\theta_P + \lambda_2 \cos\theta_P).
\end{equation}

\noindent Hence, \(\lambda_1\) and \(\lambda_2\) are constants along extremal trajectories, while \(\lambda_3\) evolves according to the equation stated above.
Characterization of the Optimal Control

\noindent The optimal control minimizes the Hamiltonian \eqref{hamiltonian} almost everywhere:
\begin{equation}\label{hmin}
    u^*(t) = \arg\min_{|u|\le a} H(x^*(t),\lambda(t),u).
\end{equation}

\noindent Since the Hamiltonian depends linearly on the control variable, define the switching function
\begin{equation}\label{switching}
    \phi(t) = \frac{\partial H}{\partial u}
= \lambda_3(t) |V_P|.
\end{equation}
\noindent Since \(|V_P|>0\), the sign of \(\phi(t)\) is determined by \(\lambda_3(t)\).
\ Then the optimal control satisfies
\begin{equation}\label{bangbang}
    u^*(t) =
\begin{cases}
-a, & \text{if } \phi(t) > 0, \\
+a, & \text{if } \phi(t) < 0.
\end{cases}
\end{equation}
\noindent If \(\phi(t)=0\) on a nontrivial interval, the control is singular and requires higher-order analysis.

\begin{theorem}
    \noindent For the interception problem defined in \eqref{problem1} and \eqref{problemT}, on every subarc on which the obstacle-avoidance constraint is inactive, optimal controls are bang–bang and take values only at the bounds \(u=\pm a\) or at \(u=0\) (straight-line motion)
\end{theorem}
\noindent Proof: A singular arc would require
\begin{equation}\label{singular}
\phi(t)=0 \quad \text{and} \quad \dot{\phi}(t)=0,
\end{equation}
\noindent over a time interval of positive measure. Since \(\phi(t)=|V_P|\lambda_3(t)\), the condition \(\dot{\phi}(t)=0\) implies
\begin{equation}\label{singular2}
-\lambda_1 \sin\theta_P + \lambda_2 \cos\theta_P = 0.
\end{equation}

\noindent Because \(\lambda_1\) and \(\lambda_2\) are constants, this identity can hold on an interval only if \(\theta_P\) is constant as \(\theta_P=\ds \tan^{-1}\frac{\lambda_2}{\lambda_1}\). From the state equation \eqref{problem1}, this requires \(u(t)=0\). Hence singular arcs correspond exclusively to straight-line segments.

\noindent Therefore, optimal trajectories consist only of:

(i). Maximum-curvature arcs (\(u=\pm a\)), and
(ii). Straight-line segments (\(u=0\)).

\section{Mathematical Model for Target Interception and Obstacle Avoidance}\label{modelM}
In this section, we present a mathematical model to determine the optimal path for a pursuer that must avoid a single obstacle—by touring around it—while intercepting a moving target. We develop mathematical expressions to compute the path lengths based on the geometric interpretation described in Section \ref{geometry1}.  The optimal interception path falls into one of four specific cases, as outlined in the set described in equation \eqref{geo_inter}, where \( CS \) denotes a segment consisting of a curve followed by a straight line, representing the pursuer’s motion, while \( S_T \) represents the straight-line path of the target.

\begin{equation}\label{pathset}
    \{LS+LS+S_T,LS+RS+S_T,RS+LS+S_T,RS+RS+S_T\}.
\end{equation}
\noindent Consider the pursuer's starting time as \( t_0 = 0 \) and the terminal time as \( t_f = t_6 \). The length of the sub-path of the pursuer is defined as \( \ell_i = |\mathbf{V_P}|(t_i - t_{i-1}) \), where \( |\mathbf{V_P}| \) is the speed of the pursuer. As the pursuer follows a path involving its minimum turning radius and tours around an obstacle, it generates six sub-paths. Therefore, we define \( \ell_i \) for \( i = 1, \ldots, 6 \). For the target, we define its path length by \( \ell_7 = |\mathbf{V_T}|(t_6 - t_0) \), where \( |\mathbf{V_T}| \) is the speed of the target, and \( \ell_7 \) represents the distance the target travels until it is intercepted by the pursuer. We solve the ordinary differential equations given in \eqref{problem1} (equivalently \eqref{kinematic}) for \( x_P(t) \), \( y_P(t) \), and \( \theta_P(t) \) over the time intervals \( t_{i-1} \leq t \leq t_i \), for \( i = 1, \ldots, 6 \).

\begin{equation}\label{solution1}
\begin{array}{cl} \ds \int_{t_{i-1}}^{t_{i}} \dot{x}_P(t) \, dt & \ =\ds \int_{t_{i-1}}^{t_{i}} |\mathbf{V_P}| \cos\theta_P (t) \,dt,  \\[4mm]
\ds \int_{t_{i-1}}^{t_{i}} \dot{y}_P(t) \, dt & \ =\ds \int_{t_{i-1}}^{t_{i}} |\mathbf{V_P}| \sin\theta_P(t) \,dt,  \\[4mm]
\ds \int_{t_{i-1}}^{t_{i}} \dot{\theta}_P(t) \, dt & \ =\ds \int_{t_{i-1}}^{t_{i}} |\mathbf{V_P}| u (t) \,dt.  \\[4mm]
\end{array}
\end{equation}

\noindent From \eqref{solution1}, we obtain the position $(x_P(t_i),y_P(t_i))$, $i=1,2$ along turning curve ($C$) yields

\begin{equation}\label{equ_6}
\begin{array}{cl} \ds x_P(t_i) & \ =\ds x_P(t_{i-1})+(\sin\theta_P(t_i)-\sin\theta_P(t_{i-1}))/(\dot{\theta}_a(t)/|\mathbf{V_P}|),  \\[4mm]
\ds y_P(t_i) & \ =\ds y_P(t_{i-1})-(\cos\theta_P(t_i)-
   \cos\theta_P(t_{i-1}))/(\dot{\theta}_a(t)/|\mathbf{V_P}|), \\[4mm]
\end{array}
\end{equation}
where 
 \begin{equation}\label{Opt.control}
\hspace{20mm}
u(t)=\frac{\dot{\theta}_a(t)}{|\mathbf{V_P}|} =
\begin{cases} 
a   &\text{if} \,\, \dot\theta_a(t) > 0 \\
-a  &\text{if}  \,\, \dot\theta_a(t) < 0 \\
0 &\text{if}\,\, \dot\theta_a(t) = 0,
\end{cases} 
\end{equation}

\noindent Assume that the centre of the fixed obstacle is $(x_b,y_b)$ and that its turning circle has radius $R_b=1/a_w$. The pursuer's entry point to the obstacle is defined as the point at which the pursuer’s straight segment meets the obstacle turning circle tangentially. This entry point is not fixed, as it depends primarily on the orientation angle $\theta_P(t)$ when the pursuer exits its turning circle. Let $\sigma\in\{+1,-1\}$ denote the touring direction: $\sigma=+1$ if the obstacle lies on the pursuer's left, so that the obstacle is toured counterclockwise by the left-turn arc $\ell_4$, and $\sigma=-1$ if the obstacle lies on the pursuer's right, so that it is toured clockwise by the right-turn arc $\ell_5$. Because the straight segment ends tangentially on the obstacle turning circle, the radius from $(x_b,y_b)$ to the entry point is orthogonal to the heading $\theta_P(t_{i-1})$, and the centre lies on the side of the pursuer given by $\sigma$. Thus, the entry coordinates $(x_P(t_i), y_P(t_i))$ for $(i=3)$ to obstacle are defined as follows. 

\begin{equation}\label{obs1}
\begin{array}{cl} \ds x_P(t_i) & \ =\ds x_b+\frac{\sigma}{a_w}\sin\theta_P(t_{i-1}), \\[4mm]
\ds y_P(t_i) & \ =\ds y_b-\frac{\sigma}{a_w}\cos\theta_P(t_{i-1}).
\end{array}
\end{equation}
\noindent Here $\sigma=\mathrm{sgn}\,w(t)$ on the touring arc, that is, $w(t)=\sigma a_w$. Hence $\sigma$ is not an additional decision variable or control; it only identifies the touring direction, which is determined by the control $w(t)$ in \eqref{Opt.control2}.

\noindent The position $(x_P(t_i),y_P(t_i)),$ with $i=3$ at which the pursuer reaches the obstacle along a straight line ($S$) can be described as

\begin{equation}\label{eq_5}
\begin{array}{cl} \ds x_P(t_i) & \ =\ds x_P(t_{i-1})+|\mathbf{V_P}|\cos\theta_P(t_{i-1})(t_i-t_{i-1}), \\[4mm]
\ds y_P(t_i) & \ =\ds y_P(t_{i-1})+|\mathbf{V_P}|\sin\theta_P(t_{i-1})(t_i-t_{i-1}). 
\end{array}
\end{equation}

\noindent Using \eqref{solution1}, we obtain the pursuer’s position \( (x_P(t_i), y_P(t_i)) \) for \( i = 4, 5 \), as it follows the turning curve \( C \) while navigating around the obstacle.

\begin{equation}\label{equ_b}
\begin{array}{cl} \ds x_P(t_i) & \ =\ds x_P(t_{i-1})+(\sin\theta_P(t_i)-\sin\theta_P(t_{i-1}))/(\dot{\theta}_b(t)/|\mathbf{V_P}|),  \\[4mm]
\ds y_P(t_i) & \ =\ds y_P(t_{i-1})-(\cos\theta_P(t_i)-
   \cos\theta_P(t_{i-1}))/(\dot{\theta}_b(t)/|\mathbf{V_P}|), \\[4mm]
\end{array}
\end{equation}
where 
 \begin{equation}\label{Opt.control2}
\hspace{20mm}
w(t)=\frac{\dot{\theta}_b(t)}{|\mathbf{V_P}|} =
\begin{cases} 
a_w   &\text{if} \,\, \dot\theta_b(t) > 0 \\
-a_w  &\text{if}  \,\, \dot\theta_b(t) < 0 \\
0 &\text{if}\,\, \dot\theta_b(t) = 0,
\end{cases} 
\end{equation}

\noindent And $i=6$ at which the pursuer intercepts the moving target along a straight line ($S$) can be described as \eqref{eq_5}


\noindent Since the speed of the pursuer ($|\mathbf{V_P}|$) is constant over the time duration $[t_{i-1}-t_{i}]$, $~~i=3, 6$, so the length of straight path for pursuer can be written as

\begin{equation}\label{solution3}
  \ds \ell_i=\ds |\mathbf{V_P}| (t_{i}-t_{i-1}), ~~i=3, 6.
\end{equation}

\noindent Now using \eqref{solution3}, we rewrite \eqref{eq_5} as below:
\begin{equation}\label{equ_7}
\begin{array}{cl} \ds x_P(t_i) & \ =\ds x_P(t_{i-1})+\ell_i\cos\theta_P(t_{i-1}), \\[4mm]
\ds y_P(t_i) & \ =\ds y_P(t_{i-1})+\ell_i\sin\theta_P(t_{i-1}), 
\end{array}
\end{equation}


\noindent and heading angles of pursuer along the both curve ($C$) and straight line ($S$) can be obtained as
\begin{equation}\label{theta_solution}
\theta_P(t_i) =
\theta_P(t_{i-1})+(\dot{\theta}_a(t)/|\mathbf{V_P}|)|\mathbf{V_P}|
  (t_i-t_{i-1}), ~~i=1,2,
\end{equation}

\noindent where $u(t)=(\dot{\theta}_a(t)/|\mathbf{V_P}|)=a$, for left-turn circular arc, $u(t)=-a$, for right-turn circular arc and $u(t)= 0$, for straight line, and

\begin{equation}\label{theta_solutionobs}
\theta_P(t_i) =
\theta_P(t_{i-1})+(\dot{\theta}_b(t)/|\mathbf{V_P}|)|\mathbf{V_P}|
  (t_i-t_{i-1}), ~~i=4,5.
\end{equation}

\noindent Similar to $u(t)$, $w(t)=(\dot{\theta}_b(t)/|\mathbf{V_P}|)$ is defined as $a_w$ for a left turn, $-a_w$ for a right turn, and $0$ for straight-line motion. Since the arc around the obstacle is centered at $(x_b,y_b)$, only the arc in the touring direction $\sigma$ is used, that is,

\begin{equation}\label{sigmacond}
\ell_5=0 \ \ \mbox{if}\ \ \sigma=+1, \qquad \ell_4=0 \ \ \mbox{if}\ \ \sigma=-1.
\end{equation}

\noindent Using \eqref{equ_6}, \eqref{Opt.control}, \eqref{obs1}, and \eqref{equ_7}(when $i=3$) the optimum path $CS$ of the pursuer for the time duration from $t_0$ to $t_3$ is obtained below.

\begin{equation}\label{path1}
    \ x_P(t_0)-\ds x_b\textcolor{blue}{-\frac{\sigma}{a_w}\sin\theta_P(t_2)}+\frac{1}{a}(-\sin\theta_P(t_0)+2\sin\theta_P(t_1)-\sin\theta_P(t_2))+\ell_3\cos\theta_P(t_2)=0,
\end{equation}
and
\begin{equation}\label{path2}
    \ y_P(t_0)-\ds y_b\textcolor{blue}{+\frac{\sigma}{a_w}\cos\theta_P(t_2)}+\frac{1}{a}(\cos\theta_P(t_0)-2\cos\theta_P(t_1)+\cos\theta_P(t_2))+\ell_3\sin\theta_P(t_2)=0.
\end{equation} 

\noindent Since the pursuer reaches the obstacle turning circle with the same heading angle so $\theta_P(t_3)=\theta_P(t_2)$. Now we generate the optimum pathway $CS$ of the pursuer when touring the obstacle for the time $t_4$ to $t_6$ using \eqref{obs1}, \eqref{equ_b}, \eqref{Opt.control2} and \eqref{equ_7}(when $i=6$).

\begin{equation}\label{path1.}
    \ds x_b\textcolor{blue}{+\frac{\sigma}{a_w}\sin\theta_P(t_2)}- x_P(t_6)+\frac{1}{a_w}(-\sin\theta_P(t_2)+2\sin\theta_P(t_4)-\sin\theta_P(t_5))+\ell_6\cos\theta_P(t_5)=0,
\end{equation}
and
\begin{equation}\label{path2.}
 \ds y_b\textcolor{blue}{-\frac{\sigma}{a_w}\cos\theta_P(t_2)}   - y_P(t_6)+\frac{1}{a_w}(\cos\theta_P(t_2)-2\cos\theta_P(t_4)+\cos\theta_P(t_5))+\ell_6\sin\theta_P(t_5)=0.
\end{equation} 
\noindent Using \eqref{problemT}, we develop the target pathway for the time duration from $t_0$  to $t_6$ as below,

\begin{equation}\label{targetint}
\int_{t_0}^{t_6}\dot{x}_T(t)\, dt=\int_{t_0}^{t_6}|\mathbf{V_T}| \cos \theta_T\, dt,
\qquad \mbox{and} \qquad
\int_{t_0}^{t_6}\dot{y}_T(t)\, dt=\int_{t_0}^{t_6}|\mathbf{V_T}| \sin \theta_T \,dt.
\end{equation}

\noindent From the above relations and using the relation $\ds \ell_7=\ds |\mathbf{V_T}| (t_6-t_0)$, we obtain the pathway for moving target as

\begin{equation}\label{path3}
    x_T(t_6)-x_T(t_0) = \ell_7\cos \theta_T,
\end{equation}
and
\begin{equation}\label{path4}
   y_T(t_6)-y_T(t_0) = \ell_7\sin \theta_T. 
\end{equation}

\noindent We combined the $CS$ pathway from obstacle to intercept point. It is noted that their positions coincide when the pursuer intercepts a moving target. In other words, $(x_P(t_6),y_P(t_6))=(x_T(t_6),y_T(t_6))$ at the interception moment. As a result combining $\eqref{path1.}$-$\eqref{path4}$ and thus we obtain,
\begin{equation}\label{path5}
    \begin{array}{l}
    x_b-x_T(t_0)\textcolor{blue}{+\ds\frac{\sigma}{a_w}\sin \theta_P(t_2)}+\frac{1}{a_w} (-\sin\theta_P(t_2)+2\sin\theta_P(t_4)-\sin\theta_P(t_5))\\[2mm]
    \qquad +\ell_6\cos\theta_P(t_5)- \ell_7 \cos\theta_T = 0,
    \end{array}
\end{equation}
and
\begin{equation}\label{path6}
    \begin{array}{l}
    y_b-y_T(t_0)\textcolor{blue}{-\ds\frac{\sigma}{a_w}\cos \theta_P (t_2)}+\frac{1}{a_w} (\cos\theta_P(t_2)-2\cos\theta_P(t_4)+\cos\theta_P(t_5))\\[2mm]
    \qquad +\ell_6\sin\theta_P(t_5)- \ell_7 \sin\theta_T = 0.
    \end{array}
\end{equation}

\noindent 
When the target is detected at time $t_0$, that point is taken as the target's initial position. At the same moment, the pursuer also begins its motion, meaning both start at the same time. For a successful interception, the pursuer and the target must reach the interception point at the same time, which means they must take the same amount of time to travel their respective paths. To ensure this timing coordination, we introduce the following constraint:

\begin{equation}\label{equ_8} \ds \sum_{i=1}^6\ell_i/ |\mathbf{V_P}|=\ell_7/|\mathbf{V_T}|.
\end{equation}
\noindent where $|\mathbf{V_T}|$ and $|\mathbf{V_P}|$ represent the constant speeds of the target and the pursuer, respectively. The terms $\ell_1$ to $\ell_6$ denote the lengths of the segments comprising the pursuer’s path, while $\ell_7$ denotes the length of the target’s path. 

\noindent Considering equations \eqref{theta_solution}-\eqref{path2} and \eqref{path5}-\eqref{equ_8}, we can formulate the following path planning model for target interception, which is the finite-dimensional counterpart of the optimal control problem \eqref{problem1} for the path structure in Figure \ref{Figopti_path}.

\begin{equation}\label{model_1}
 \begin{array}{cl} \min &\ f =\ds \sum_{i=1}^7 \ell_i \\[2mm]
 \mbox{subject to} \\[2mm]
 & \ x_P(t_0)-\ds x_b-\frac{\sigma}{a_w}\sin\theta_P(t_2)+\frac{1}{a}(-\sin\theta_P(t_0)+2\sin\theta_P(t_1)-\sin\theta_P(t_2)) \\[2mm] & +\ell_3\cos\theta_P(t_2)=0,\\[2mm]
 & \ y_P(t_0)-\ds y_b+\frac{\sigma}{a_w}\cos\theta_P(t_2)+\frac{1}{a}(\cos\theta_P(t_0)-2\cos\theta_P(t_1)+\cos\theta_P(t_2))  \\[2mm] & +\ell_3\sin\theta_P(t_2)=0,\\ [2mm]
 & \ x_b-x_T(t_0)+\ds\frac{\sigma}{a_w}\sin \theta_P(t_2)+\frac{1}{a_w} (-\sin\theta_P(t_2)+2\sin\theta_P(t_4)-\sin\theta_P(t_5)) \\[2mm] & +\ell_6\cos\theta_P(t_5)- \ell_7 \cos\theta_T = 0,\\ [2mm]
 & \ y_b-y_T(t_0)-\ds \frac{\sigma}{a_w}\cos \theta_P (t_2)+\frac{1}{a_w} (\cos\theta_P(t_2)-2\cos\theta_P(t_4)+\cos\theta_P(t_5)) \\[2mm] & +\ell_6\sin\theta_P(t_5)- \ell_7 \sin\theta_T = 0, \\ [2mm]
 & \ds |\mathbf{V_T}|\sum_{i=1}^6\ell_i - |\mathbf{V_P}|\ell_7 =0,\\ [2mm]
&\ \ell_i \geq 0, \;\;  i=1,\ldots,7, \quad \ell_1,\ell_2\leq 2\pi/a, \quad \ell_4,\ell_5\leq 2\pi/a_w,
\end{array} 
\end{equation}
where 
\begin{equation}\label{model1theta}
\begin{array}{ll}
    \theta_{P}(t_1)=\theta_{P}(t_0) + a \ell_1, & \theta_{P}(t_2) =\theta_{P}(t_1) - a \ell_2, \\[1mm]
    \theta_{P}(t_4)=\theta_{P}(t_2)  + a_w \ell_4, & \theta_{P}(t_5)  =\theta_{P}(t_4)  - a_w \ell_5,
\end{array}
\end{equation}
and $\sigma$ and the pair $(\ell_4,\ell_5)$ satisfy \eqref{sigmacond}. The upper bounds on the arc lengths exclude arcs of more than one full turn, which would only repeat the same configuration. Each of the four path types in \eqref{pathset} is obtained by fixing the initial turn and the touring direction: $\ell_2=0$ for an initial left turn ($L$) and $\ell_1=0$ for an initial right turn ($R$), and $\sigma=+1$ ($\ell_5=0$) for touring with a left turn and $\sigma=-1$ ($\ell_4=0$) for touring with a right turn. Model \eqref{model_1} is solved for each of the four path types, and the feasible path of minimum length is selected as the optimal path. Because $\ds \sum_{i=1}^{6}\ell_i=(|\mathbf{V_P}|/|\mathbf{V_T}|)\,\ell_7$ by the time-synchronisation constraint, minimizing $\ds f=\sum_{i=1}^{7}\ell_i$ is equivalent to minimizing the pursuer's path length \eqref{minilengthpath} and the interception time $\ds t_6=\sum_{i=1}^{6}\ell_i/|\mathbf{V_P}|$.

\begin{remark}\label{rem:inactive}
The tangency conditions in \eqref{model_1} replace the obstacle-avoidance inequality of \eqref{problem1} by the geometric structure of its boundary arc: the pursuer reaches the obstacle turning circle tangentially, follows it with curvature $\pm a_w$, and leaves it tangentially. Two cases are therefore checked when the model is applied. (i) If the obstacle-free interception path (an initial arc followed by a straight segment to the interception point, obtained from \eqref{model_1} without the obstacle terms and with $\ell_4=\ell_5=\ell_6=0$) does not intersect the obstacle turning circle, the obstacle-avoidance constraint is inactive and this path is optimal; the obstacle need not be toured. (ii) Otherwise, the candidate paths of the four types in \eqref{pathset} are computed, any candidate whose arcs or straight segments enter the obstacle turning circle is discarded, and the shortest remaining candidate is selected.
\end{remark}

\subsection{Touring $n$ static obstacles}
\noindent We extend Model \eqref{model_1} to the case where the pursuer must avoid $n$ static obstacles before intercepting a moving target. The goal is to minimize the total path length, allowing the pursuer to intercept the target along the shortest possible path. The path is divided into segments made up of $CS$ (curve-straight) type paths. The first $CS$ path goes from the pursuer’s starting point to the entry point of the first obstacle. The second $CS$ path connects the entry point of the first obstacle to the entry point of the second obstacle, and this process continues for each subsequent obstacle. In total, $ n$ $ CS$ paths are required to reach the entry point of the $n$th obstacle. Another path $CS$ is added from the entry point of the last obstacle to the point of intercept with the target. This final segment requires initial target position and direction information. So, with $n$ obstacles, the entire path includes $(n+1)$ $CS$ paths, which are divided into $3n + 4$ segments, labeled $\ell_i$, where $i = 1, \dots, 3n + 4$. As in the single-obstacle case, each obstacle is entered and left tangentially, and the collision checks of Remark \ref{rem:inactive} are applied to every obstacle.

\noindent This extends the model for touring $n$ obstacles and intercepting a moving target as follows.

\begin{equation}\label{model_2}
 \begin{array}{cl} \min &\ f =\ds \sum_{i=1}^{3n+4} l_i, \;\;\; n\in \mathbb{N}\\[2mm]
 \mbox{subject to} & \\ &
   \begin{cases}
  & \ x_P(t_{3j-3})-\ds x_P(t_{3j})+\frac{1}{a_j}(-\sin\theta_P(t_{3j-3})+2\sin\theta_P(t_{3j-2})-\sin\theta_P(t_{3j-1})) \\ & + \;\; \ell_{3j}\cos\theta_P(t_{3j})=0,\\[2mm]
 & \ y_P(t_{3j-3})-\ds y_P(t_{3j})+\frac{1}{a_j}(\cos\theta_P(t_{3j-3})-2\cos\theta_P(t_{3j-2}) +\cos\theta_P(t_{3j-1})) \\ & + \;\; \ell_{3j}\sin\theta_P(t_{3j})=0,
\end{cases}\\[5mm]
 & \ x_P(t_{3n})-x_T(t_0)+\ds \frac{1}{a_{n+1}} \left(-\sin\theta_P(t_{3n})+2\sin\theta_P(t_{3n+1}\right)-\sin\theta_P(t_{3n+2}))\\
 & +\ell_{3n+3}\cos\theta_P(t_{3n+3})- \ell_{3n+4} \cos\theta_{T} = 0,\\ [2mm]
 & \ y_P(t_{3n})-y_T(t_0)+\ds\frac{1}{a_{n+1}} \left(\cos\theta_P(t_{3n})-2\cos\theta_P(t_{3n+1})+\cos\theta_P(t_{3n+2})\right)\\
 & +\ell_{3n+3}\sin\theta_P(t_{3n+3})- \ell_{3n+4} \sin\theta_{T} = 0, \\ [2mm]
 & \ |\mathbf{V_T}| \ds \sum_{i=1}^{3n+3}\ell_i - |\mathbf{V_P}| \ell_{3n+4}=0,\;\; \ell_i \geq 0,  \; \mbox{for} \,\; i=1,\ldots ,3n+4.
\end{array} 
\end{equation}
where 

\begin{equation}\label{model2entry}
    x_P(t_{3j})=x_c^j+\ds \frac{\sigma_j}{a_{j+1}}\sin\theta_P(t_{3j}), \;\; \mbox{and}\;\; y_P(t_{3j})=y_c^j-\ds \frac{\sigma_j}{a_{j+1}}\cos\theta_P(t_{3j}),\;\;j=1,\ldots,n,
\end{equation}
where $\sigma_j\in\{+1,-1\}$ is the touring direction of the $j$th obstacle, with $\ell_{3j+2}=0$ if $\sigma_j=+1$ and $\ell_{3j+1}=0$ if $\sigma_j=-1$,
\begin{equation}\label{model2theta}
    \theta_P(t_j)=\theta_{P}(t_{j-1})+\alpha_j l_j, \;\; j=1, \ldots, 3n+3. 
\end{equation}
\begin{equation}\label{model2alpha}
\hspace{-40mm}
    \begin{cases}
  \ds \;\;\alpha_{3k-2}=\,\,\,\,a_k,\\[2mm]
    \ds \;\;\alpha_{3k-1}=-a_k,\;\; \\[2mm]
  \ds \;\;  \alpha_{3k}=\,\,\,\,0,\;\;\\[2mm]
   \end{cases}
\end{equation}
where,   $\ds \;\; k \in \{1, 2, \ldots, n+1\}$.

\section{ Numerical Experiments}\label{NM}
The results of our numerical experiments are presented in this section to verify the effectiveness of our proposed model \eqref{model_2} in determining the optimal route that a pursuer must take to intercept a moving target interacting with static obstacles. To determine the best route for target interception, we provide the initial positions and heading angles of the pursuer and the target, as well as the obstacle's center and turning radii, taking into account their constant speeds.  In each experiment, the pursuer has feasible routes, and our model determines the optimum pathway among the paths. We also evaluate our proposed model by varying the positions, speeds, and turning radii of both the pursuer and the target. This allows us to assess the efficiency of the model and examine how well it reaches the objectives of the analysis.

\noindent We perform the experiments using AMPL \cite{Fourer2003} modelling software, employing a variety of solvers, including GUROBI, SCIP, LGO, BONMIN, and IPOPT \cite{Wachter2006}, all with their default settings. MATLAB is used to simulate the solutions obtained from these numerical experiments. All simulations are carried out on an HP ProBook $440$ $G6$ laptop equipped with a $4.6$ GHz Intel Core $i7$ processor and $8$ GB of RAM.

\noindent \textbf{Example 1}\label{example_1}

\noindent In this example, three experimental scenarios are carried out to evaluate the model’s effectiveness in computing optimal and feasible trajectories, as shown in Table \ref{table:experiment-1a}.
  \par
\begin{figure}[H]
\hspace{-1cm}
\begin{minipage}{100mm}
\begin{center}
\hspace{-1.5cm}
\includegraphics[width=90mm]{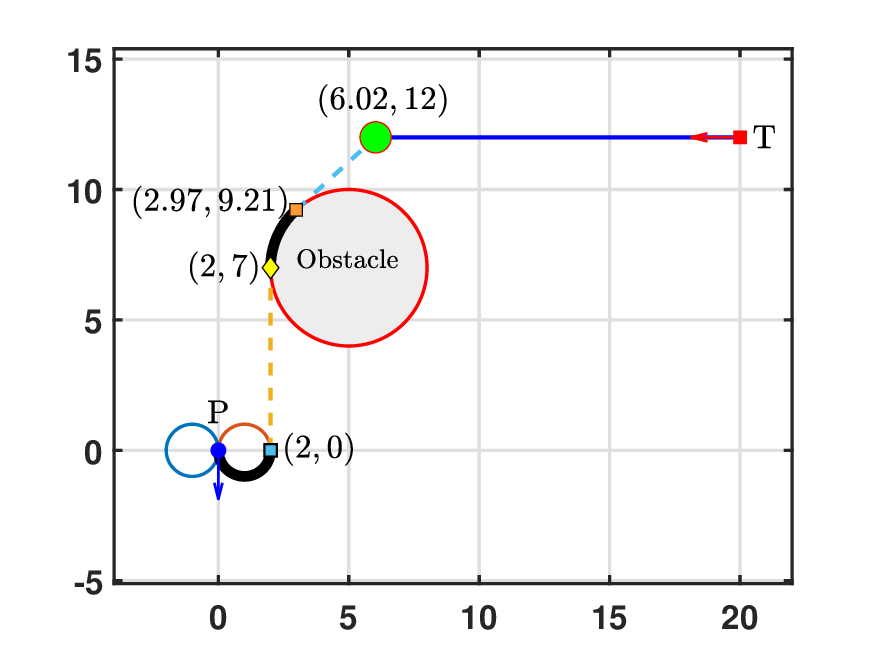} \\
\vspace*{-0.3cm}
{\scriptsize (a) Optimal path.}
\end{center}
\end{minipage}
\hspace{-2.6cm}
\begin{minipage}{100mm}
\begin{center}
\hspace*{0cm}
\includegraphics[width=90mm]{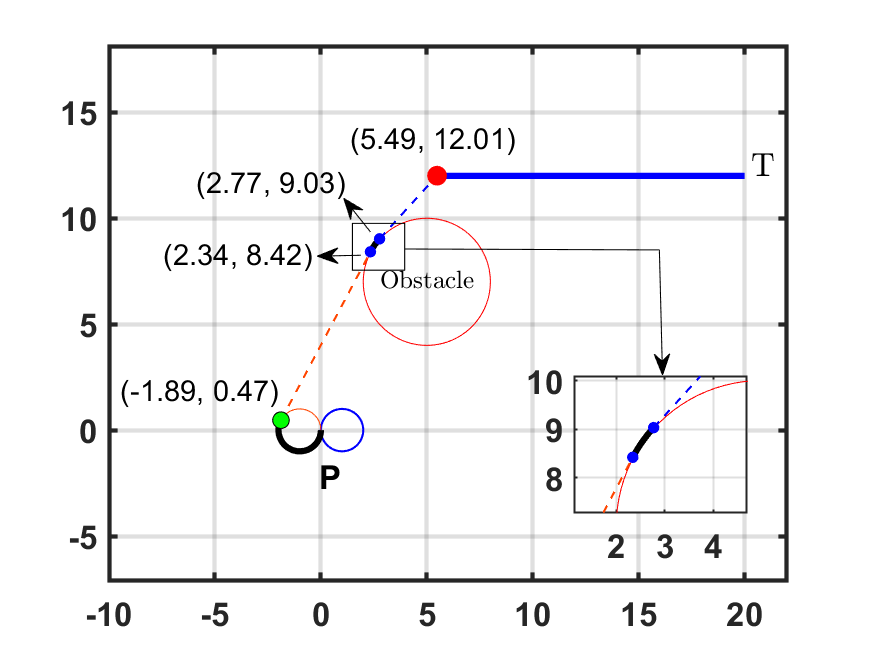} \\
\vspace*{-0.3cm}
{\scriptsize (b) Feasible path.}
\end{center}
\end{minipage}
\vspace*{.4cm}
\caption{\sf\footnotesize For Table \ref{table:experiment-1a} (a), optimal and feasible paths illustrating the exit points of the initial turning circles, the entry and exit points of subsequent turning circles around obstacles, and the interception points of the target for both paths.}
\label{EX1}
\end{figure}

	\begin{table}[H]
		\caption{\small{\textit{Test cases with a single obstacle. }}}
		\footnotesize
		\vskip 1.5em
		\centering
		\begin{tabular}{|c| c| c| c| c| }
  \hline
\multicolumn{1}{|c|}{} &\multicolumn{1}{|c|}{Pursuer}& \multicolumn{1}{|c|}{Target}  & \multicolumn{1}{|c|}{Obstacle} & \multicolumn{1}{|c|}{} \\ & & & (centre) &\\
\hline
			\hline
           & $(0,0,3\pi/2)$ & $(20,12,\pi)$ & $(5,7)$ & Position\\ 
           $(a)$ & 6 & 5 & & Speed\\
           & 1 & &3 & Turning Radius\\
            [.5ex]\hline
          &  $(0,0,7\pi/6)$ & $(30,15,\pi)$ & $(6,8)$ & Position\\ 
           $(b)$& 6 & 2 & & Speed\\
           & 0.5 & &5 & Turning Radius\\
            [.5ex]\hline
          &  $(0,0,7\pi/6)$ & $(6,8,\pi/4)$ & $(6,8)$ & Position\\ 
          $(c)$ &  6 & 3 & & Speed\\
           & 1 & &5 & Turning Radius\\
            [.5ex]\hline
            
		\end{tabular}
		\label{table:experiment-1a}
	\end{table}
\noindent In Table \ref{table:experiment-1a}(a), we assume that the pursuer initiates its trajectory from the position \( P(0,0, 3\pi/2) \), while the target's initial position is set at \( T(20,12, \pi) \). A static obstacle is centered at \( (5,7) \) with a turning radius of $3$. In this experiment, the pursuer maintains a constant speed of $6$ units, whereas the target moved along a straight path at a speed of $5$ units. This analysis demonstrated that the optimal trajectory for the pursuer follows the sequence \( LS+RS+S_T \) (see Figure \ref{EX1}(a)): the pursuer enters the obstacle turning circle tangentially at $(2,7)$, tours it clockwise, and leaves it at $(2.97, 9.21)$. Figure \ref{EX1}(b) presents the second-best feasible path, of type \( RS+RS+S_T \), corresponding to the parameter settings provided in Table \ref{table:experiment-1a}(a).

\noindent Table \ref{table:exp-1}(a) shows the optimal path taken by the pursuer to intercept the target. The path is composed of several segments, including left turns, right turns, and straight-line movements. The pursuer starts with a left turn of length $\bar{\ell}_1 = 3.14$, followed by a straight segment of $\bar{\ell}_3 = 7.00$. Next, the pursuer makes a right turn of $\bar{\ell}_5 = 2.49$ around the obstacle, then continues along a straight segment of $\bar{\ell}_6 = 4.13$ to reach the target. Meanwhile, the target moves along a straight path of $\bar{\ell}_7 = 13.97$ before being intercepted. The total length of the combined path, including both the pursuer and target segments, is $f = 30.74$, representing the shortest path that satisfies the motion constraints and interception conditions. Detailed information about the feasible path is provided in Table \ref{table:exp-1}(a) under the row labeled ``Feasible Path".

  \par
	\begin{table}[H]
		\caption{\small{\textit{Segment lengths of optimal and feasible paths for the setup outlined in Table \ref{table:experiment-1a}. }}}
		\footnotesize
		\vskip 1.5em
		\centering
		\begin{tabular}{|c| c| c| c| c| c| c| c|c| c| }
  \hline
\multicolumn{1}{|c|}{}& \multicolumn{1}{|c|}{}& \multicolumn{7}{|c|}{Segment Lengths}  & \multicolumn{1}{|c|}{} \\
\hline
			\hline
			& & Left  &  Right  & St.	& Left  &  Right  & St.  & Target  & Total \\
			& &     turn		&   turn  	& line  & turn& turn &line & St. line & length   \\
			& &$\bar{\ell}_1$ &     $\bar{\ell}_2$ 		&  $\bar{\ell}_3$   	& $\bar{\ell}_4$ & $\bar{\ell}_5$& $\bar{\ell}_6$ & $\bar{\ell}_7$ & $f$ \\ [0.5ex]
			\hline
			$(a)$ & Optimal Path & 3.14  & 0 & 7.00 & 0   &	2.49  & 4.13& 13.97 	& 30.74  \\ [.5ex]\hline
            & Feasible Path & 0  & 3.63 & 9.00 & 0   &	0.75  & 4.03& 14.51 	& 31.91  \\ [.5ex]\hline
           $(b)$ & Optimal Path & 1.57  & 0 & 9.16 & 1.18   &	0  & 15.62& 9.17 	& 36.69  \\ [.5ex]\hline
           & Feasible Path & 5.10 & 0 & 13.93  &  0  &	7.25  & 8.35 & 11.54 & 46.17  \\ [.5ex]\hline
           $(c)$ & Optimal Path & 0  & 2.41 & 8.78 & 0   &	4.14  & 13.54& 14.43 	& 43.30  \\ [.5ex]\hline
						
		\end{tabular}
		\label{table:exp-1}
	\end{table}

\noindent Table \ref{table3:exp-1} presents the time taken for each path section, including left turns, right turns, and straight-line segments, for both the optimal and feasible paths. Based on the parameters provided in Table \ref{table:experiment-1a}(a), Table \ref{table3:exp-1}(a) shows that the optimal path takes in total \textcolor{blue}{$2.79$} seconds to intercept the target. It is noted that the computational time to determine the optimal route in this experiment was approximately $0.031$ seconds.

  \par
\begin{figure}[H]
\hspace{-1cm}
\begin{minipage}{100mm}
\begin{center}
\hspace{-1.5cm}
\includegraphics[width=95mm]{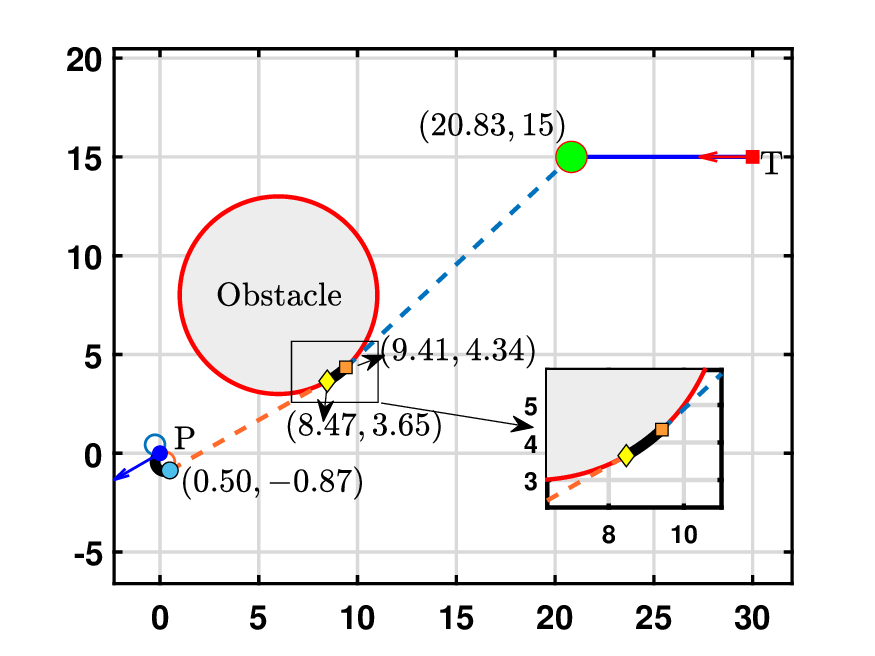} \\
\vspace*{-.3cm}
{\scriptsize (a) Optimal path.}
\end{center}
\end{minipage}
\hspace{-2.0cm}
\begin{minipage}{100mm}
\begin{center}
\hspace*{0cm}
\includegraphics[width=95mm]{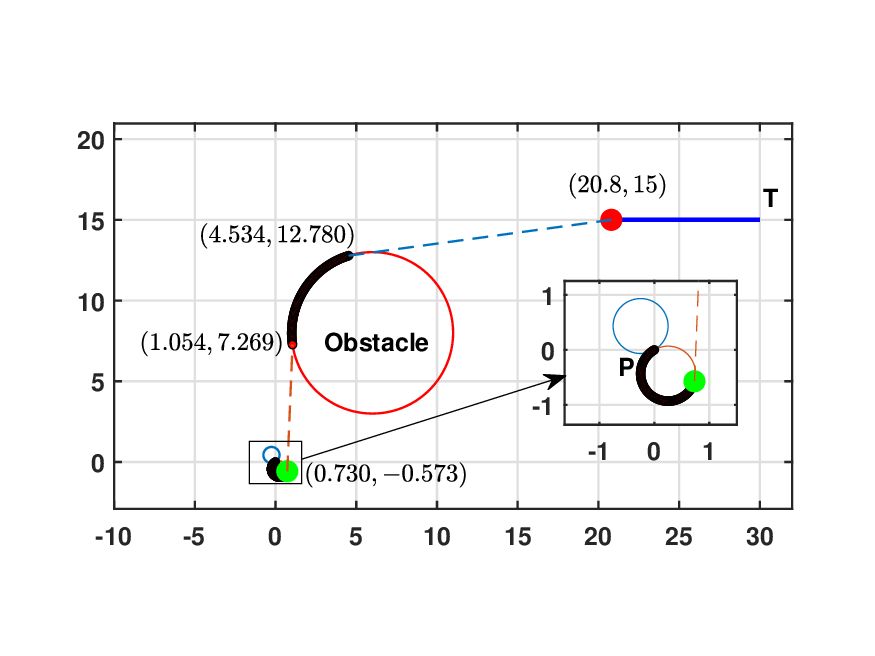} \\
\vspace*{-.9cm}
{\scriptsize (b) Feasible path.}
\end{center}
\end{minipage}
\vspace*{.4cm}
\caption{\sf\footnotesize For Table \ref{table:experiment-1a} (b), optimal and feasible paths illustrating the exit points of the initial turning circles, the entry and exit points of subsequent turning circles around an obstacle, and the interception points of the target for both paths.}
\label{EX2}
\end{figure}

	\begin{table}[!htbp]
		\caption{\small{\textit{Path segment durations based on Table \ref{table:experiment-1a}. }}}
		\footnotesize
		\vskip 1.5em
		\centering
		\begin{tabular}{|c| c| c| c| c| c| c| c|c| c| }
  \hline
\multicolumn{1}{|c|}{} & \multicolumn{1}{|c|}{}& \multicolumn{7}{|c|}{Time}  & \multicolumn{1}{|c|}{} \\
\hline
			\hline
			& & Left  &  Right  & St.	& Left  &  Right  & St.  & Target  & Total \\
           &  &     turn		&   turn  	& line  & turn& turn &line & St. line & Time   \\
			& &$\bar{t}_1$ &     $\bar{t}_2$ 		&  $\bar{t}_3$   	& $\bar{t}_4$ & $\bar{t}_5$& $\bar{t}_6$ & $\bar{t}_7$ & $t$ \\ [0.5ex]
			\hline
			$(a)$ & Optimal Time & 0.52  & 0 & 1.17 & 0   &	0.42  & 0.69& 2.79 	& 2.79  \\ [.5ex]\hline
            & Feasible Time & 0  & 0.61 & 1.50 & 0   &	0.12  & 0.67& 2.90 	& 2.90  \\ [.5ex]\hline
           $(b)$ & Optimal Time & 0.26  & 0 & 1.53 & 0.20   &	0  & 2.60& 4.59 	& 4.59  \\ [.5ex]\hline
           & Feasible Time & 0.85  & 0 & 2.32  &  0  &	1.21  & 1.39 & 5.77  	& 5.77  \\ [.5ex]\hline
          $(c)$ &  Optimal Time & 0  & 0.40 & 1.46 & 0   &	0.69 & 2.26& 4.81 	& 4.81  \\ [.5ex]\hline
						
		\end{tabular}
		\label{table3:exp-1}
	\end{table}
\noindent We now test our model \eqref{model_2} using the parameter settings given in Table \ref{table:experiment-1a} (b) and (c), where the pursuer’s starting positions and heading angles are varied along with the target and obstacle locations. For the configuration in Table \ref{table:experiment-1a}(b), the analysis shows that the optimal trajectory for the pursuer follows the sequence $LS + LS + S_T$, as illustrated in Figure \ref{EX2}(a), while the corresponding feasible path of type $LS+RS+S_T$ is shown in Figure \ref{EX2}(b). The total optimal path length is $f = 36.69$, as reported in Table \ref{table:exp-1}(b), and the total optimal time taken is $t = 4.59$ seconds, shown in Table \ref{table3:exp-1}(b). The arc around the obstacle is short in this case ($\bar{\ell}_4=1.18$); it is enlarged in the inset of Figure \ref{EX2}(a).

\begin{figure}[H]
\vspace*{-15mm}
 \centering
\includegraphics[width=100mm]{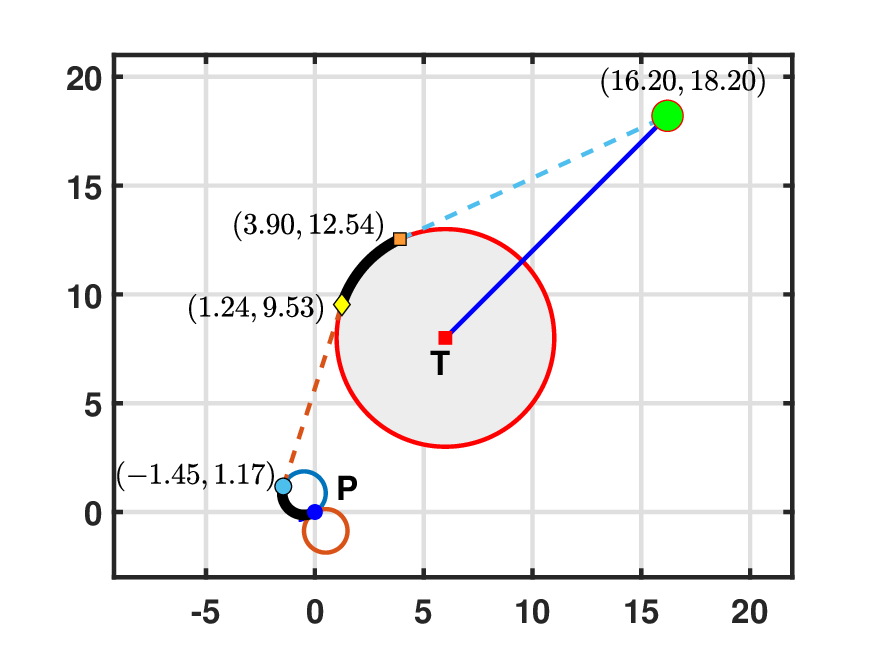} \\ \vspace{-.5cm}
{\footnotesize  Optimal path}.
\\
 \vspace{.4cm}
\caption{\sf\footnotesize For Table \ref{table:experiment-1a} (c), optimal path illustrating the exit points of the initial turning circles, the entry and exit points of subsequent turning circles around an obstacle, and the interception points of the target for optimal path.}
\label{EX3}
\end{figure}

\noindent For the parameter settings in Table \ref{table:experiment-1a}(c), the results are illustrated in Figure \ref{EX3}, with corresponding path length and time details provided in Tables \ref{table:exp-1}(c) and \ref{table3:exp-1}(c). The optimal path is of type $RS+RS+S_T$, with $f=43.30$ and interception time $t=4.81$ seconds.
\newpage





\noindent \textbf{Example 2} \label{example2}

\noindent We test the proposed Model \eqref{model_2} to handle cases involving multiple obstacles. In this experiment, we generate optimal paths in the presence of two obstacles across four different setups, as listed in Table \ref{table:experiment-2a}.

  \par
	\begin{table}[H]
		\caption{\small{\textit{Test cases with two obstacles. }}}
		\footnotesize
		\vskip 1.5em
		\centering
		\begin{tabular}{|c| c| c| c| c| c|}
  \hline
\multicolumn{1}{|c|}{}& \multicolumn{1}{|c|}{Pursuer}& \multicolumn{1}{|c|}{Target}  & \multicolumn{1}{|c|}{Obstacle 1} & \multicolumn{1}{|c|}{Obstacle 2} &  \multicolumn{1}{|c|}{} \\ & & & (centre) & (centre) &\\
\hline
			\hline
          &  $(0,0,\pi)$ & $(0,25,0)$ & $(6,8)$ & $(15,15)$ & Position\\ 
          $(a)$ &  7 & 5 & & & Speed\\
           & 1 & &4 &4 & Turning Radius\\
            [.5ex]\hline
          &  $(0,0,7\pi/6)$ & $(25,5,\pi/2)$ & $(6,8)$ & $(15,20)$ & Position\\ 
         $(b)$ &   6 & 3 & & & Speed\\
           & 1 & &5 &4 & Turning Radius\\
            [.5ex]\hline
          &  $(0,0,\pi)$ & $(25,25,3\pi/2)$ & $(6,8)$ & $(15,15)$ & Position\\ 
          $(c)$ &  7 & 2 & & & Speed\\
          &  1 & &4 &4 & Turning Radius\\
            [.5ex]\hline
           &  $(0,0,5\pi/6)$ & $(35,5,4\pi/6)$ & $(6,8)$ & $(15,15)$ & position\\ 
         $(d)$ &   7 & 2 & & & Speed\\
           & 1 & &4 &4 & Turning Radius\\
            [.5ex]\hline
            
		\end{tabular}
		\label{table:experiment-2a}
	\end{table}

\begin{figure}[!htbp]
\centering
\begin{minipage}[b]{0.45\textwidth}
\centering
\includegraphics[width=\textwidth,trim=205 30 196 34,clip]{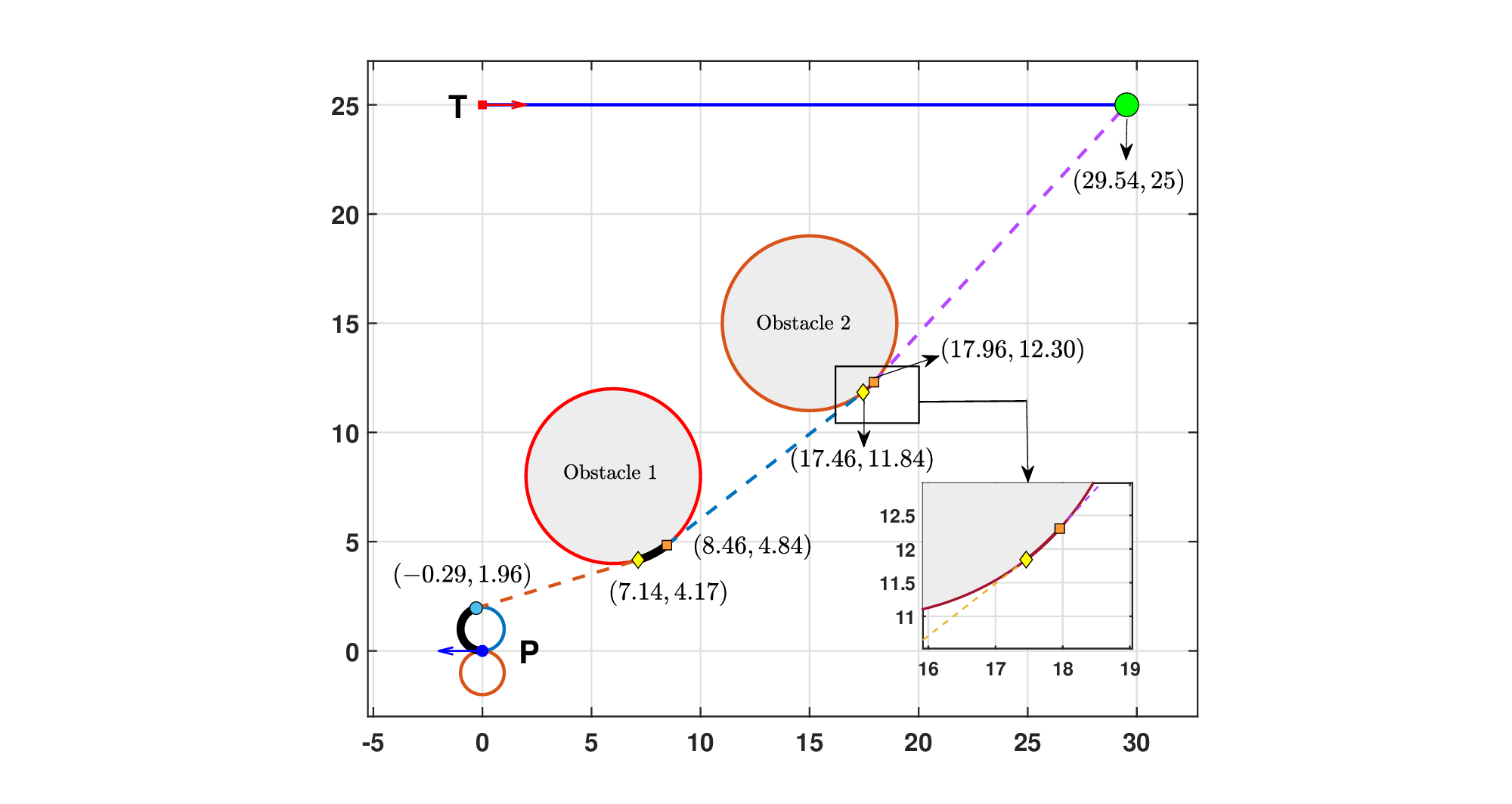}\\[1mm]
{\scriptsize (a) Optimal path for Table \ref{table:experiment-2a}(a).}
\end{minipage}\hfill
\begin{minipage}[b]{0.45\textwidth}
\centering
\hspace{-1.1cm}
\includegraphics[width=\textwidth,trim=205 30 196 34,clip]{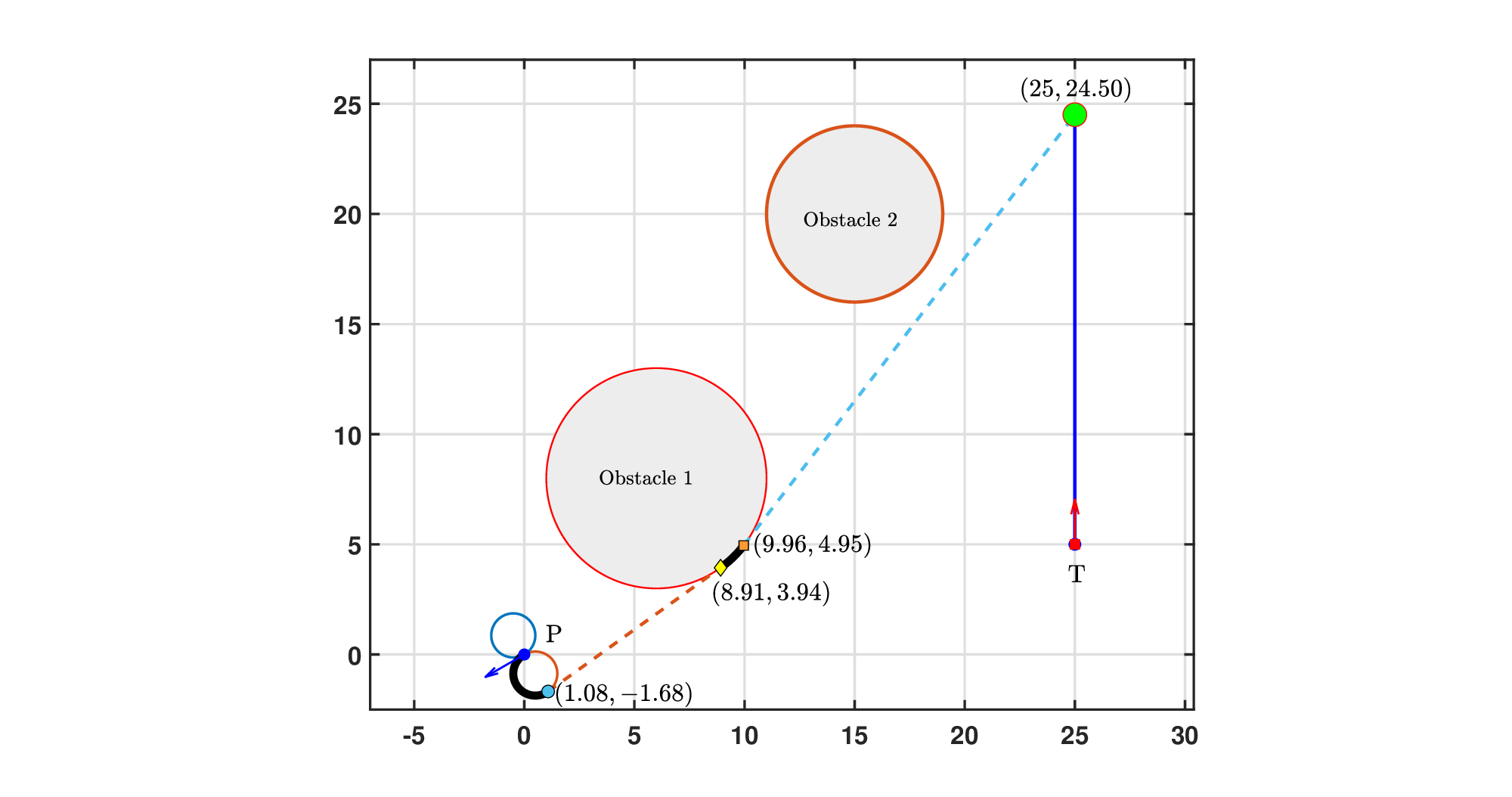}\\[1mm]
{\scriptsize (b) Optimal path for Table \ref{table:experiment-2a}(b).}
\end{minipage}\\[4mm]
\begin{minipage}[b]{0.49\textwidth}
\centering
\includegraphics[width=\textwidth,trim=205 30 196 34,clip]{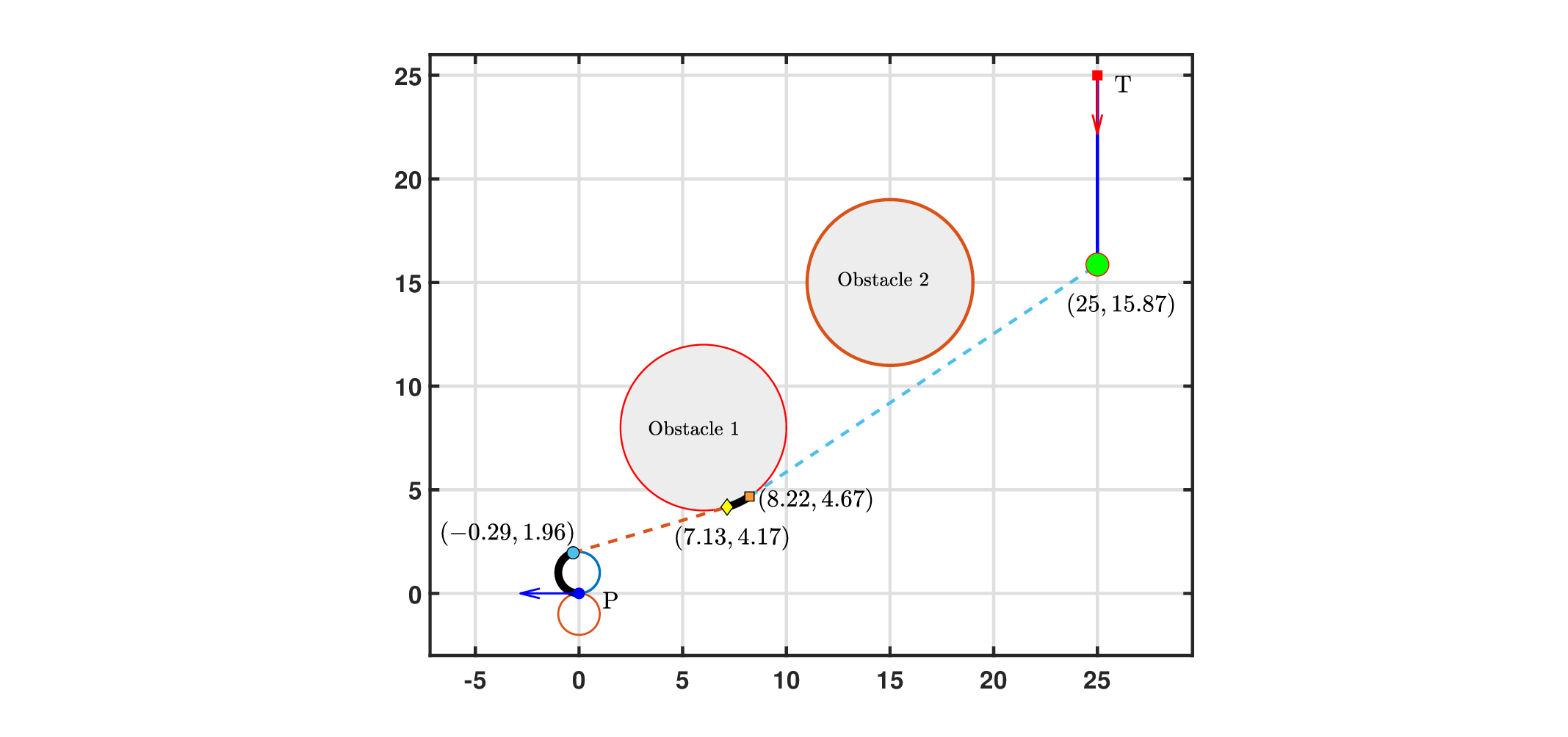}\\[1mm]
{\scriptsize (c) Optimal path for Table \ref{table:experiment-2a}(c).}
\end{minipage}\hfill
\begin{minipage}[b]{0.49\textwidth}
\centering
\includegraphics[width=\textwidth,trim=205 30 196 34,clip]{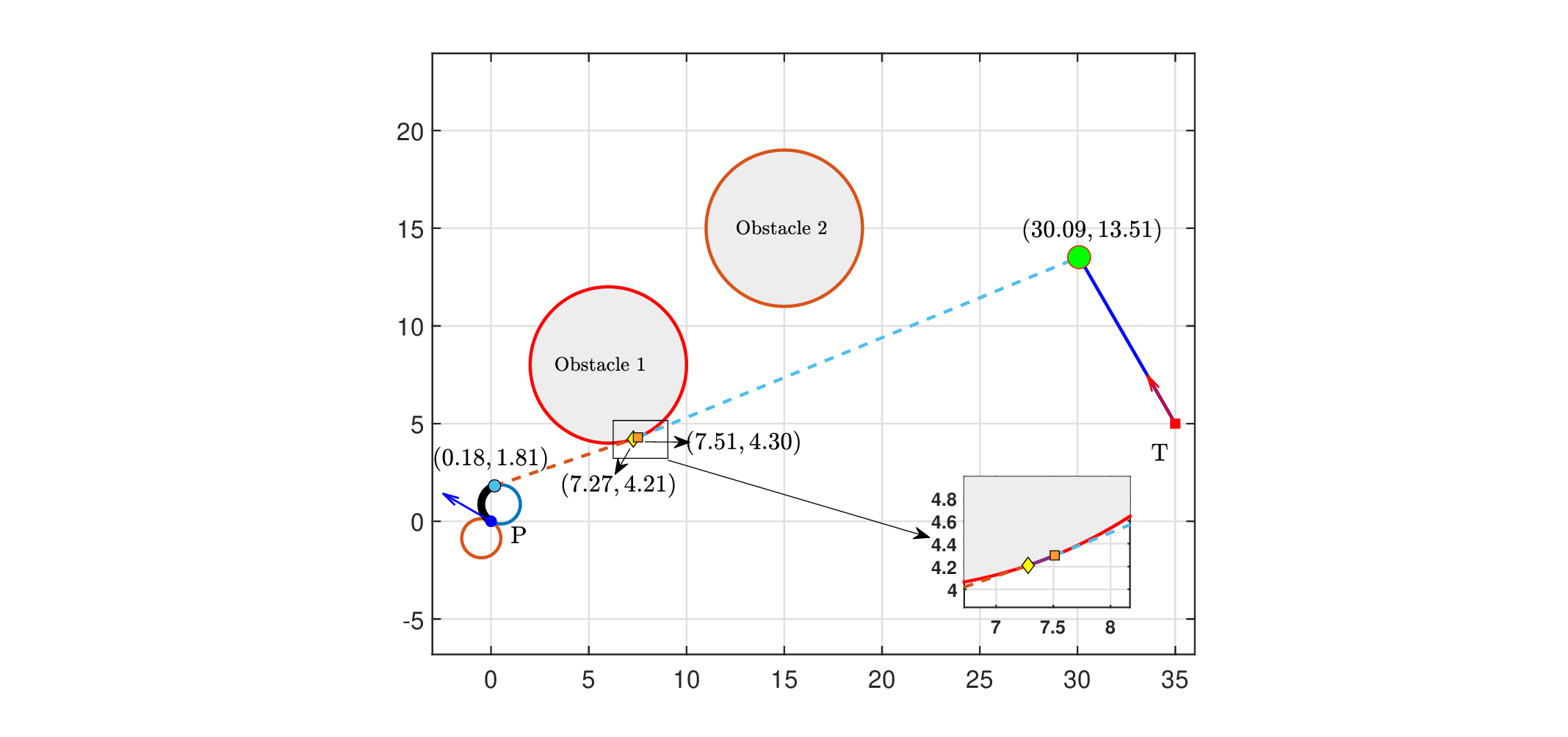}\\[1mm]
{\scriptsize (d) Optimal path for Table \ref{table:experiment-2a}(d).}
\end{minipage}
\vspace*{2mm}
\caption{\sf\footnotesize For Table \ref{table:experiment-2a}(a)-(d), optimal paths illustrating the exit points of the initial turning circles, the entry and exit points of subsequent turning circles around two obstacles, and the interception points of the target.}
\label{EX2b}
\end{figure}

\noindent For Table \ref{table:experiment-2a}(a), the pursuer's trajectory originates from $P(0,0, \pi)$, while the target begins its motion from $T(0,25, 0)$. Two static obstacles, centered at $(6,8)$ and $(15,15)$, each with a turning radius of $4$, are also incorporated into the analysis. The pursuer maintains a constant speed of $7$ units, whereas the target moves along a straight path at a speed of $5$ units. Our analysis determines that the optimal trajectory for the pursuer follows the sequence $RS+LS+LS+S_T$ demonstrated in Figure \ref{EX2b}(a). 


\noindent 
Table \ref{table:exp-4}(a) presents the optimal path taken by the pursuer to intercept the target. The path consists of a sequence of segments—left turns, right turns, and straight-line movements similar to the structure observed in Example 1 with a single obstacle. The pursuer begins with a right turn of length $\bar{\ell}_2 = 2.85$, followed by a straight segment of $\bar{\ell}_3 = 7.75$. It then makes a left turn of $\bar{\ell}_4 = 1.49$ around the first obstacle, continues with a straight segment of $\bar{\ell}_6 = 11.40$ toward the second obstacle, makes a short left turn of $\bar{\ell}_7 = 0.68$ around it (enlarged in the inset of Figure \ref{EX2b}(a)), and follows a straight segment of $\bar{\ell}_9 = 17.19$ to the interception point. Meanwhile, the target moves along a straight path of $\bar{\ell}_{10} = 29.54$ before being intercepted. The total combined path length, including both the pursuer and target segments, is $f = 70.90$, representing the shortest optimal path that meets the motion constraints and interception requirements.

  \par
	\begin{table}[H]
		\caption{\small{\textit{Segment lengths of optimal paths for the setup outlined in Table \ref{table:experiment-2a}. }}}
		\footnotesize
		\vskip 1.5em
		\centering
		\begin{tabular}{|c| c| c| c| c| c| c| c| c| c| c|c| c| }
  \hline
\multicolumn{1}{|c|}{}& \multicolumn{1}{|c|}{}& \multicolumn{10}{|c|}{Length}  & \multicolumn{1}{|c|}{} \\
\hline
			\hline
			& & Left  &  Right  & St.	& Left  &  Right  & St. & Left  &  Right  & St. & Target  & Total \\
         &  &     turn		&   turn  	& line  & turn& turn &line & turn& turn &line & St. line & length   \\
		&	&$\bar{\ell}_1$ &     $\bar{\ell}_2$ 		&  $\bar{\ell}_3$   	& $\bar{\ell}_4$ & $\bar{\ell}_5$& $\bar{\ell}_6$& $\bar{\ell}_7$ & $\bar{\ell}_8$& $\bar{\ell}_9$ & $\bar{\ell}_{10}$ & $f$ \\ [0.5ex]
			\hline
		$(a)$	& Optimal Path & 0  & 2.85 & 7.75 & 1.49   &	0  & 11.40& 0.68 & 0& 17.19 &29.54 	& 70.90  \\ [.5ex]\hline
        $(b)$  &  Optimal Path & 3.24  & 0 & 9.64 & 1.47   &	0  & 24.67& 0 & 0& 0 &19.50	& 58.51  \\ [.5ex]\hline
         $(c)$ &  Optimal Path & 0  & 2.85 & 7.75 & 1.20   &	0  & 20.17& 0 & 0& 0 &9.13 	& 41.10  \\ [.5ex]\hline
        $(d)$  &  Optimal Path & 0  & 2.29 & 7.49 & 0.25   &	0  & 24.38& 0 & 0& 0 & 	9.83& 44.25 \\ [.5ex]\hline
          &  Feasible Path 1 & 4.92  & 0 & 13.67 & 0   &	5.75  & 1.83& 1.45 & 0& 15.98 &12.60 	& 56.20 \\ [.5ex]\hline
          &  Feasible Path 2 & 0  & 1.48 & 11.39 & 0   &	2.23  & 7.07& 0 & 3.11& 14.90 &11.48 	& 51.66 \\ [.5ex]\hline
						
		\end{tabular}
		\label{table:exp-4}
\\[1mm] \parbox{0.95\textwidth}{\footnotesize When a path does not tour the second obstacle, $\bar{\ell}_7=\bar{\ell}_8=\bar{\ell}_9=0$ and $\bar{\ell}_6$ is the final straight segment to the interception point.}
	\end{table}

\noindent Other setups are presented in Table \ref{table:experiment-2a}(b), (c), and (d), where the pursuer and target heading angles, target positions, and the location of the second obstacle are varied. Using these configurations, we test the proposed Model \eqref{model_2}. The corresponding optimal paths for each case are illustrated in Figures \ref{EX2b}(b), (c), and (d), where the entry and exit points on the turning circles are indicated. Detailed sub-length information for these optimal paths is provided in Table \ref{table:exp-4}(b), (c), and (d), respectively. A summary of the simulation results—including interception points, total path lengths, and time required for interception—for the cases listed in Table \ref{table:experiment-2a} is presented in Table \ref{table:experiment-2b}. In cases (b)--(d), the optimal path tours only the first obstacle: the second obstacle does not obstruct the remaining straight segment to the interception point, so its avoidance constraint is inactive (Remark \ref{rem:inactive}), although the path is still checked against it.

  \par
	\begin{table}[H]
		\caption{\small{\textit{Simulation result summary of Table \ref{table:experiment-2a}. }}}
		\footnotesize
		\vskip 1.5em
		\centering
		\begin{tabular}{|c| c| c| c| c|}
  \hline
\multicolumn{1}{|c|}{}&\multicolumn{1}{|c|}{Optimum Path}& \multicolumn{1}{|c|}{Interception Point}  & \multicolumn{1}{|c|}{Optimum Length} & \multicolumn{1}{|c|}{Optimum Time} \\
\hline
			\hline
         $(a)$ &  $RS+LS+LS+S_T$ & $(29.54,25)$ & $f=70.90$ & $t=5.91$\\[.5ex]\hline
         $(b)$ &  $LS+LS+S_T$ & $(25,24.50)$ & $f=58.51$ & $t=6.50$\\[.5ex]\hline
         $(c)$ &  $RS+LS+S_T$& $(25,15.87)$ & $f=41.10$ & $t=4.57$\\[.5ex]\hline
         $(d)$ & $RS+LS+S_T$& $(30.09,13.51)$ & $f=44.25$ & $t=4.92$\\[.5ex]\hline
                  
		\end{tabular}
		\label{table:experiment-2b}
	\end{table}

\begin{figure}[H]
\hspace{-1cm}
\begin{minipage}{100mm}
\begin{center}
\hspace{-1.5cm}
\includegraphics[width=90mm]{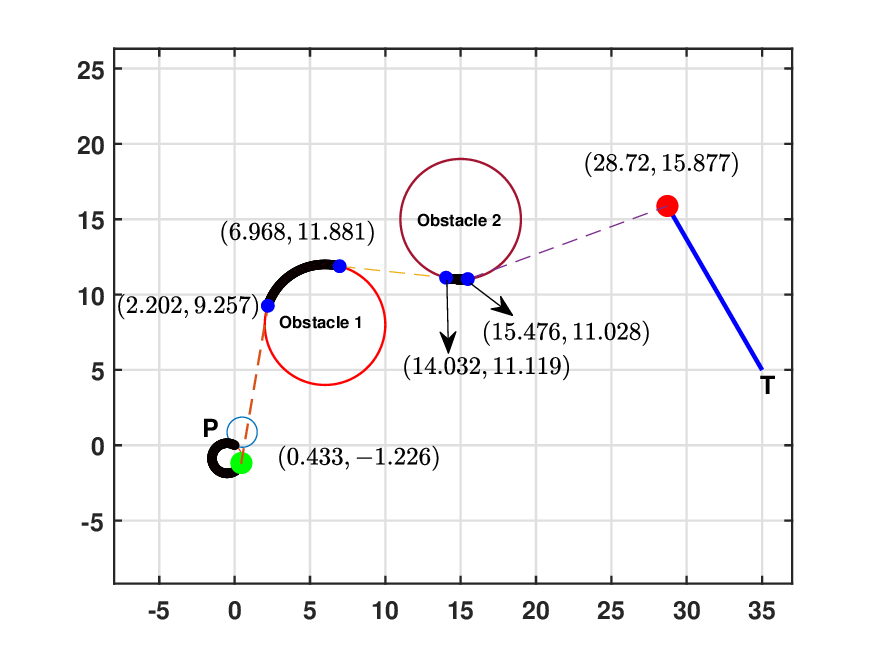} \\
\vspace*{-.5cm}
{\scriptsize (a) Feasible path 1.}
\end{center}
\end{minipage}
\hspace{-2.6cm}
\begin{minipage}{100mm}
\begin{center}
\hspace*{0cm}
\includegraphics[width=90mm]{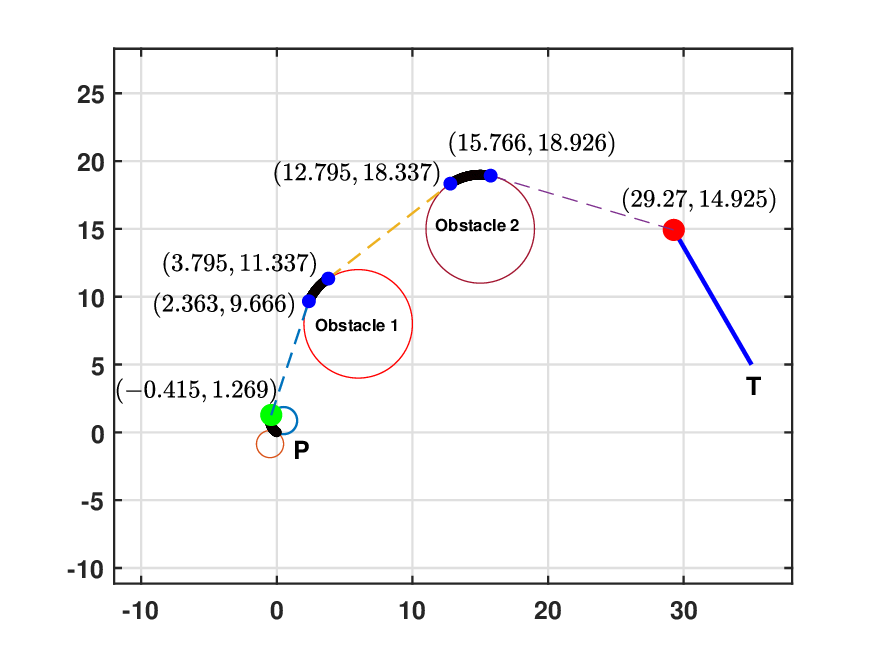} \\
\vspace*{-.5cm}
{\scriptsize (b) Feasible path 2.}
\end{center}
\end{minipage}
\vspace*{.4cm}
\caption{\sf\footnotesize For Table \ref{table:experiment-2a}(d), feasible paths illustrating the exit points of the initial turning circles, the entry and exit points of subsequent turning circles around two obstacles, and the interception points of the target for both paths.}
\label{EX2c}
\end{figure}

\noindent Two feasible paths are presented in Figure \ref{EX2c} for the parameter setting in Table \ref{table:experiment-2a}(d), alongside the optimal path shown in Figure \ref{EX2b}(d). These paths are included to illustrate that while the optimal path minimizes total length, alternative trajectories still satisfy the model’s constraints. The first feasible path (Figure \ref{EX2c}(a), $f=56.20$) is the path $LS+RS+LS+S_T$, which passes both obstacles, and the second (Figure \ref{EX2c}(b), $f=51.66$) is of type $RS+RS+RS+S_T$ and tours both obstacles. Both are longer than the optimal path $RS+LS+S_T$ with $f=44.25$ (Table \ref{table:exp-4}(d)). Showing feasible paths helps demonstrate the model’s flexibility and the impact of different trajectory choices.
\noindent Examples $1$ and $2$, along with supporting experiments, demonstrate the model’s ability to compute optimal paths under varying initial conditions while successfully avoiding obstacles and intercepting the target. The model effectively determines the shortest interception time and also provides feasible alternative trajectories.


\noindent \textbf{Example 3 (Challenging scenarios)}\label{example3}

\begin{figure}[!htbp]
\centering
\begin{minipage}[b]{0.48\textwidth}
\centering
\includegraphics[width=\textwidth,trim=219 32 187 14,clip]{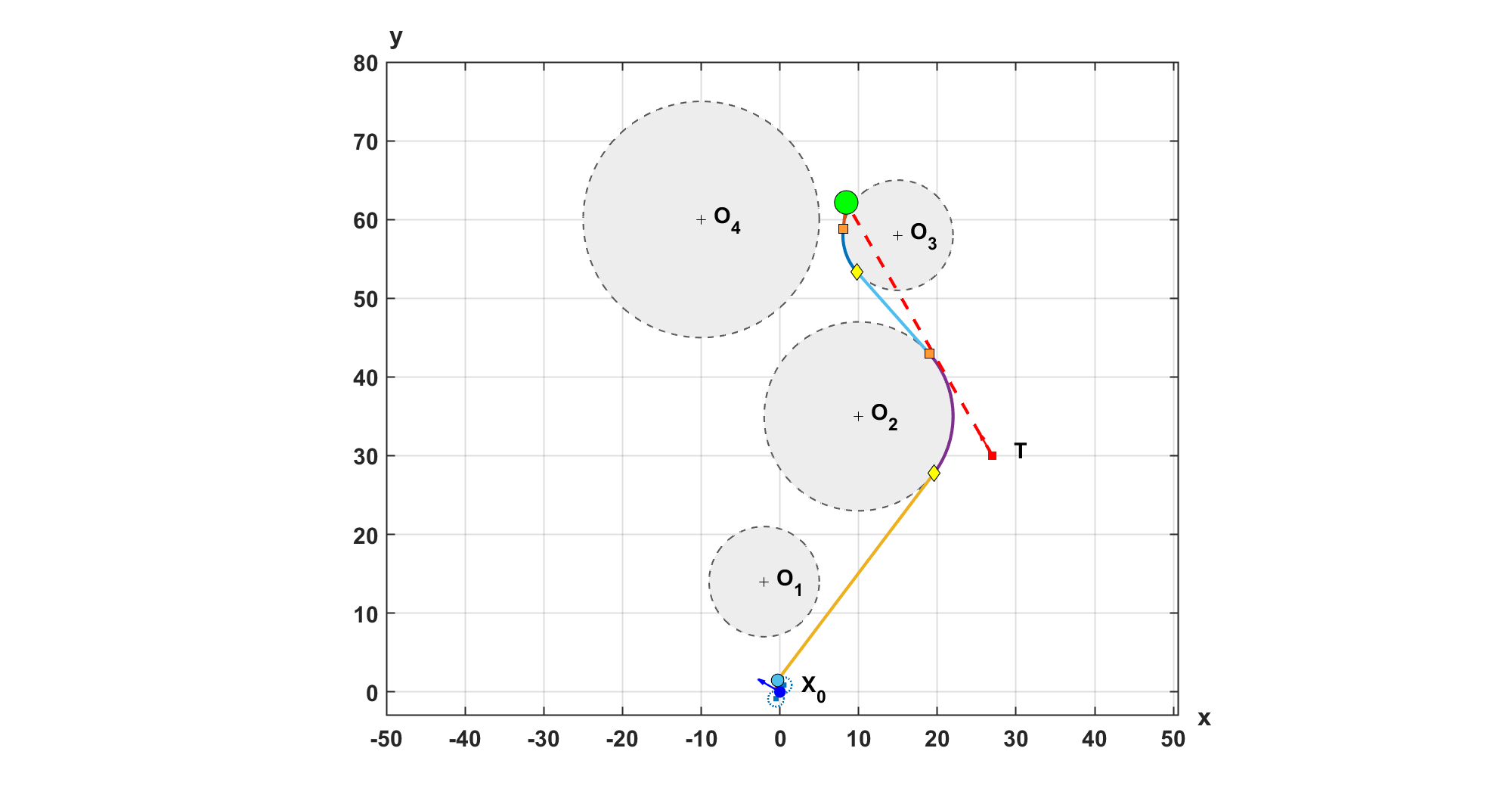}\\[1mm]
{\footnotesize (a) $|\mathbf{V_T}|=3$.}
\end{minipage}\hfill
\begin{minipage}[b]{0.48\textwidth}
\centering
\includegraphics[width=\textwidth,trim=140 27 269 18,clip]{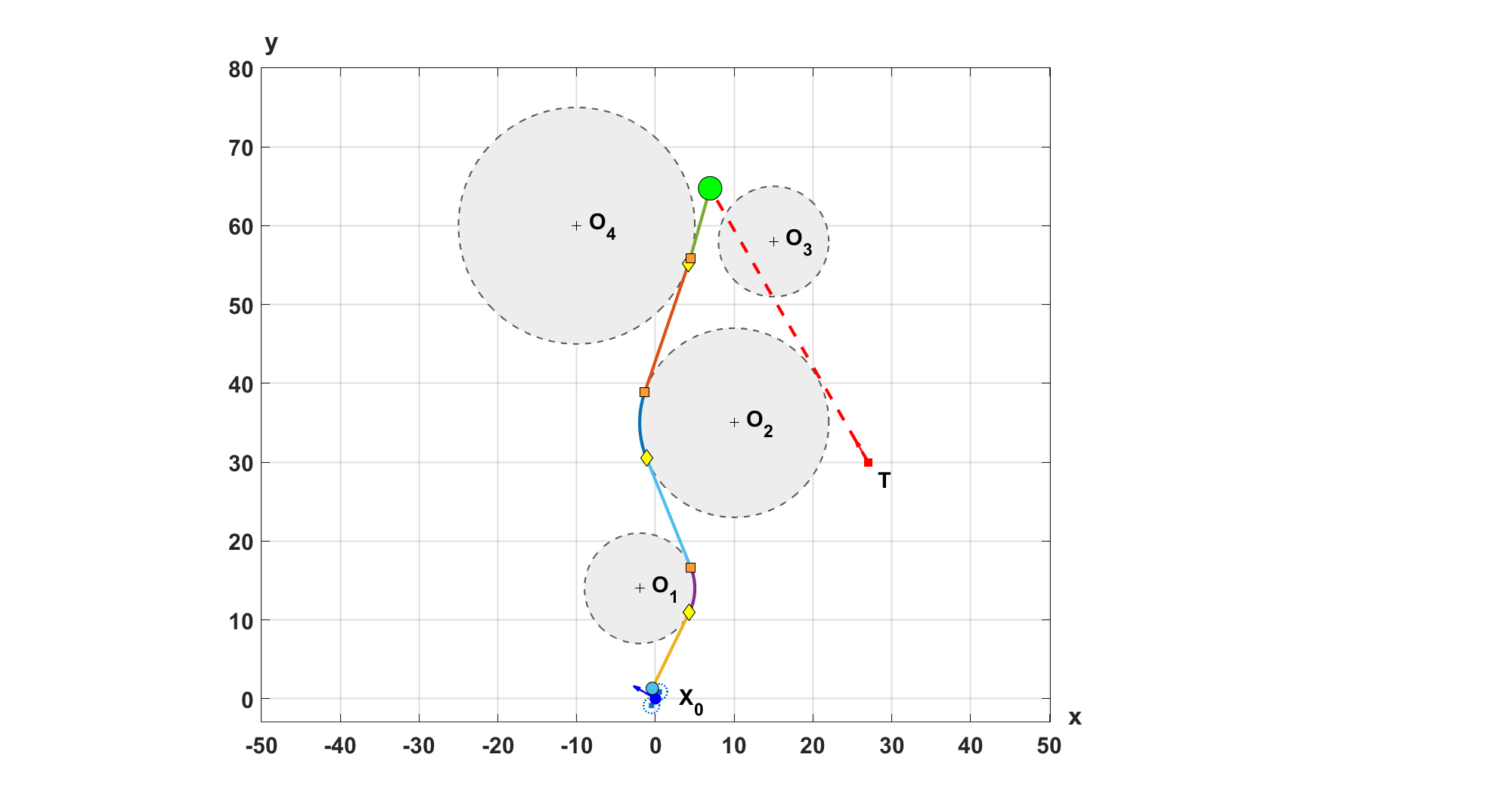}\\[1mm]
{\footnotesize (b) $|\mathbf{V_T}|=3.5$.}
\end{minipage}\\[4mm]
\begin{minipage}[b]{0.48\textwidth}
\centering
\includegraphics[width=\textwidth,trim=207 31 185 16,clip]{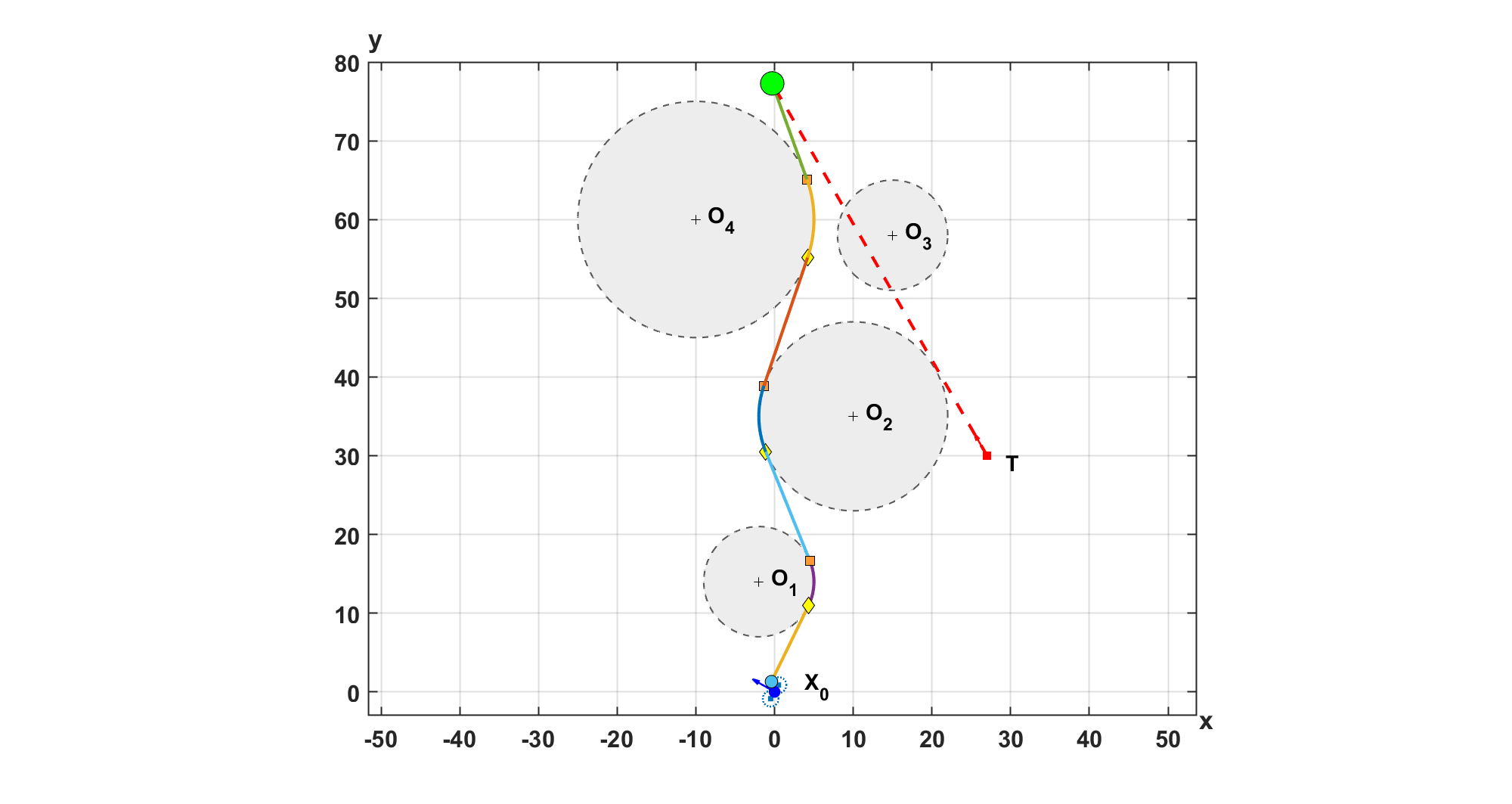}\\[1mm]
{\footnotesize (c) $|\mathbf{V_T}|=4$.}
\end{minipage}
\vspace*{2mm}
\caption{\sf\footnotesize Optimal interception paths among four obstacles $O_1$--$O_4$ for target speeds (a) $3$, (b) $3.5$ and (c) $4$, with all other parameters fixed. The markers indicate the initial point $X_0$ of the initial turning circle, the entry (diamond) and exit (square) points on the obstacle turning circles, and the interception point (green circle).}
\label{fig:speeds}
\end{figure}

\noindent Under the revised model formulation, we further tested the model and the solution algorithm on several challenging scenarios. Here the pursuer must pass through narrow gaps between four static obstacles of different sizes, $O_1$--$O_4$, before it can intercept the target. The pursuer starts from $P(0,0)$ and the target from $T$, and all parameters except the target speed are kept fixed. Figure \ref{fig:speeds} shows the optimal paths for target speeds $|\mathbf{V_T}|=3$, $3.5$ and $4$; for brevity, these results are presented only graphically.

\noindent In all three cases, the pursuer enters and leaves each obstacle turning circle tangentially, and no part of the path enters an obstacle. The obstacles toured, and the side on which each one is passed, change with the target speed. For $|\mathbf{V_T}|=3$, the pursuer passes $O_1$ without touring it and goes around the right-hand side of $O_2$. For $|\mathbf{V_T}|=3.5$ and $4$, it tours $O_1$ and then passes $O_2$ on its left-hand side before traveling through the gap between $O_3$ and $O_4$. As the target speed increases, the interception point moves further along the target's path. These results show that the model, together with the enumeration of path types and the collision check of Remark \ref{rem:inactive}, selects a feasible shortest path even when several obstacles lie close to each other and the best route around them changes with the target's motion.

\subsection{Comparison with existing interception and obstacle-avoidance methods}\label{sec:comparison}

\noindent Table \ref{tab:qualitative_comparison} compares the proposed formulation with representative Dubins-based approaches for intercepting a moving target. Meyer et al. \cite{Meyer2015} and Zheng et al. \cite{Zheng2021} study the interception of a target moving at constant velocity by a Dubins vehicle with constant speed and a minimum turning radius, and characterize the minimum-length (equivalently, minimum-time) interception paths analytically, but they do not consider obstacles. Triharminto et al. \cite{Triharminto2014} consider a UAV intercepting a moving target in the presence of an obstacle and propose a turning strategy based on a modified Dubins algorithm, which is evaluated by simulation without an optimality guarantee.

\begin{table}[H]
\centering
\caption{\small{\textit{Qualitative comparison of the proposed model with existing Dubins-based interception methods.}}}
\label{tab:qualitative_comparison}
\footnotesize
\renewcommand{\arraystretch}{1.2}
\begin{tabular}{|>{\raggedright\arraybackslash}p{3.3cm}|>{\raggedright\arraybackslash}p{2.3cm}|>{\raggedright\arraybackslash}p{2.3cm}|>{\raggedright\arraybackslash}p{2.5cm}|>{\raggedright\arraybackslash}p{2.8cm}|}
\hline
\textbf{Feature} & \textbf{Meyer et al. \cite{Meyer2015}} & \textbf{Zheng et al. \cite{Zheng2021}} & \textbf{Triharminto et al. \cite{Triharminto2014}} & \textbf{Proposed model} \\
\hline
Moving target (constant velocity) & $\checkmark$ & $\checkmark$ & $\checkmark$ & $\checkmark$ \\ \hline
Constant speed and minimum turning radius & $\checkmark$ & $\checkmark$ & $\checkmark$ & $\checkmark$ \\ \hline
Terminal heading & free & free & not specified & free \\ \hline
Static obstacles & $\times$ & $\times$ & $\checkmark$ (single) & $\checkmark$ (finite number) \\ \hline
Obstacle handling & -- & -- & modified Dubins turning rule & tangential touring of the obstacle turning circle \\ \hline
Interception / time synchronisation & analytical & analytical & simulation-based & explicit constraints \eqref{path5}--\eqref{equ_8} \\ \hline
Solution approach & geometric / analytical & PMP-based guidance law & rule-based planner & nonlinear programme via arc parameterisation \\ \hline
Optimality criterion & minimum length / time & minimum time & not stated & minimum length (= minimum time) \\ \hline
\end{tabular}
\end{table}

\noindent The comparison shows that the proposed formulation combines moving-target interception and the avoidance of a finite number of static obstacles within one optimization model, in which the obstacle-touring geometry, the interception point and the time synchronisation are represented explicitly. When the obstacle does not obstruct the interception path, the obstacle-avoidance constraint is inactive (Remark \ref{rem:inactive}(i)) and the model reduces to the obstacle-free interception problem studied in \cite{Meyer2015,Zheng2021}. Because the methods in Table \ref{tab:qualitative_comparison} address different problem settings (with or without obstacles, analytical guidance laws or rule-based planning), the comparison is intended to highlight the additional problem components handled by the proposed model rather than to rank the methods numerically. A quantitative benchmark against sampling-based planners adapted to moving-target interception is left for future work.

\section{Real-World Test Scenario}\label{RWA}

\noindent In this test example, we utilise location data sourced from Google Maps, including the locations and heading angles of both the pursuer and the target, as well as the positions and radii (sizes) of obstacles that the pursuer must avoid. The experiment is set in Bangladesh using Google Maps. After defining the locations of all relevant elements, we collect the necessary data to carry out the simulation. The UAV (pursuer) starts from Satkhira City to intercept a moving target that begins traveling in a straight line from Feni City. Two restricted airspaces (or circular obstacles) are defined over Narail and Narayanganj, which the pursuer must strictly avoid. The latitude and longitude coordinates for all relevant locations are collected using Google Maps. Initial parameters are summarized in Table \ref{table:Application1}.  The autonomous pursuer is modeled with a minimum turning radius of \(10\) km.  We used a relatively large minimum turning radius because the speed was set to 800 km/hour. However, we also tested the algorithm with smaller turning radii to evaluate its performance under tighter manoeuvring conditions. Each restricted airspace (obstacle) is represented with a turning radius of \(20\) km, defining the area that the pursuer must avoid.

  \par
	\begin{table}[!htbp]
		\caption{\small{\textit{Locations, coordinates, and motion parameters }}}
		\footnotesize
		\vskip 1.5em
		\centering
		\begin{tabular}{|c| c| c| c| c| }
  \hline
\multicolumn{1}{|c|}{Name of City}& \multicolumn{1}{|c|}{Position in Earth}  & \multicolumn{1}{|c|}{Position in 2D} & \multicolumn{1}{|c|}{Speed} &  \multicolumn{1}{|c|}{Turning Radius} \\ & & & km/hour & \\
\hline
			\hline
            Satkhira & $
22^\circ 44' 15.66'' \text{N}, \quad 89^\circ 3' 25.65'' \text{E}$
 & $(9133.10,2528.32)$ & $800$ & 10\\[.5ex]\hline 
            Narail& $
23^\circ 11' 44.15'' \text{N}, \quad 89^\circ 29' 49.47'' \text{E}
$ &$(9147.19,2579.23)$ & & 20\\[.5ex]\hline
           Narayanganj &$
23^\circ 38' 35.46'' \text{N}, \quad 90^\circ 28' 55.86'' \text{E}
$ &$(9216.63,2629.00)$ & & 20\\
            [.5ex]\hline
            Feni & $
23^\circ 1' 34.54'' \text{N}, \quad 91^\circ 22' 43.76'' \text{E}
$&$(9351.30,2560.40)$ & $650$&\\[.5ex]\hline
            
		\end{tabular}
		\label{table:Application1}
	\end{table}

\begin{figure}[H]
\hspace{-1cm}
\begin{minipage}{100mm}
\begin{center}
\hspace{-1.5cm}
\includegraphics[width=85mm]{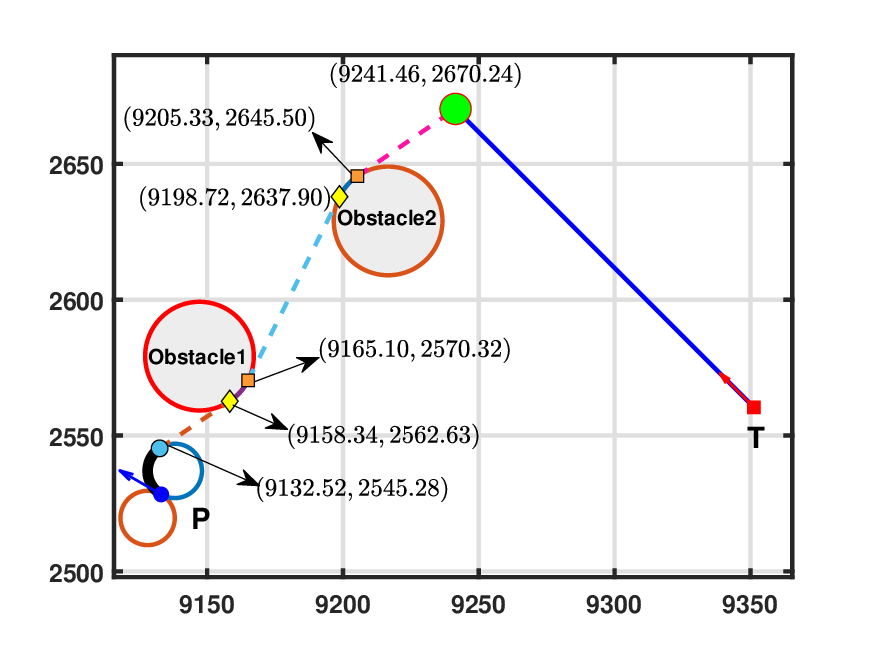} \\
\vspace*{-.3cm}
{\scriptsize (a) Optimal path.}
\end{center}
\end{minipage}
\hspace{-2.6cm}
\begin{minipage}{100mm}
\begin{center}
\hspace*{-1cm}
\includegraphics[width=80mm]{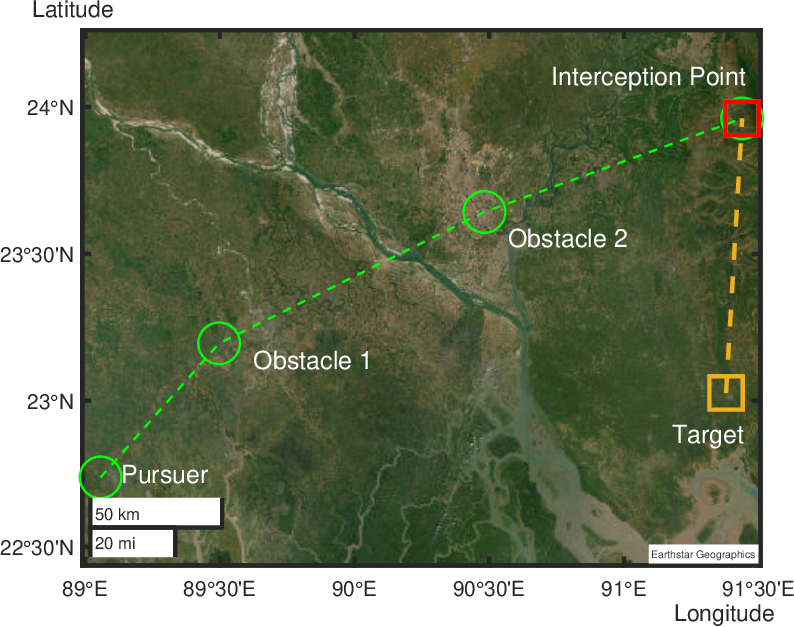} \\
\vspace*{-0.3cm}
{\scriptsize (b) Geographic Locations.}
\end{center}
\end{minipage}
\vspace*{.4cm}
\caption{\sf\footnotesize (a) Optimized trajectory of the automated drone from P$(9133.10,2528.32)$ with heading $\theta_P = 5\pi/6$, intercepting the target that moved from $T(9351.30,2560.40)$ while avoiding the obstacles centered at $(9147.19,2579.23)$ and $(9216.63,2629.00)$  as the target drone moves in direction $\theta_T = 3\pi/4$. (b) Geographic Locations on Map.}
\label{App1}
\end{figure}

\noindent  In our analysis, the target moves with heading $\theta_T=3\pi/4$ and the pursuer successfully intercepts the target at geographic coordinates
approximately $24^\circ 0' 50.6'' \, \text{N}$, $90^\circ 58' 20.3'' \, \text{E}$,
which correspond to 2D coordinates $(9241.46, 2670.24)$. A detailed breakdown of the optimal interception path length is provided in Table \ref{table:application2}. The optimal path is of type $RS+LS+RS+S_T$: the pursuer tours the restricted airspace over Narail counterclockwise and the one over Narayanganj clockwise. The pursuer covers approximately $191.18$ km to reach the target, while the target travels about $155.33$ km before being intercepted, so that the interception occurs after $0.239$ hours (about $14.3$ minutes). The total interception distance, $f = 346.51 \, \text{km}$, is the sum of both trajectories. The simulation scenario is illustrated in Figure \ref{App1}, which also depicts the corresponding geographical setting.
.

  \par
	\begin{table}[H]
		\caption{\small{\textit{Optimized path generation under model performance. }}}
		\footnotesize
		\vskip 1.5em
		\centering
		\begin{tabular}{|c| c| c| c| c| c| c| c| c| c| c|c| }
  \hline
\multicolumn{1}{|c|}{$\mathbf{|V|}$ $\&$ $R$}& \multicolumn{10}{|c|}{Length}  & \multicolumn{1}{|c|}{} \\
\hline
			\hline
			$\mathbf{|V_P|}=800$ & Left  &  Right  & St.	& Left  &  Right  & St. & Left  &  Right  & St. & Target  & Total \\
			$\mathbf{|V_T|}=650$ &     turn		&   turn  	& line  & turn& turn &line & turn& turn &line & St. line & length   \\
			$\;~R\;\;=10$&$\bar{\ell}_1$ &     $\bar{\ell}_2$ 		&  $\bar{\ell}_3$   	& $\bar{\ell}_4$ & $\bar{\ell}_5$& $\bar{\ell}_6$& $\bar{\ell}_7$ & $\bar{\ell}_8$& $\bar{\ell}_9$ & $\bar{\ell}_{10}$ & $f$ \\ [0.5ex]
			\hline
			Optimal Path & 0  & 20.26 & 31.11 & 10.35   &	0  & 75.49& 0 & 10.18& 43.78 &155.33 	& 346.51  \\ [.5ex]\hline

		\end{tabular}
		\label{table:application2}
	\end{table}



\section{Conclusion}
\noindent We have presented an optimal path planning approach for intercepting a moving target while avoiding obstacles. The study focuses on a $2$D optimal control problem for a pursuer, considering both kinematic and tactical constraints. A geometric configuration is introduced to illustrate how the sub paths are determined. The pursuer's path is subject to a constraint that ensures the turning radius remains greater than or equal to its minimum turning capability. The mathematical formulation begins with the case of a single obstacle and is later extended to accommodate a finite number of obstacles. We derived necessary optimality conditions, and the model is illustrated through extensive numerical experiments. Several illustrative simulations are presented, demonstrating the interception path for various combinations of initial positions and heading angles of both the pursuer and the target. Additionally, a real-world case study is included to validate that the proposed optimal path planning method is capable of handling interception tasks in environments with obstacles.

\noindent \textbf{Limitations and future work:} Despite the promising theoretical and numerical results, the proposed approach has several limitations. First, the model assumes constant speeds for both the pursuer and the target, with the target moving along a known straight-line trajectory. Second, the obstacle-avoidance formulation assumes static obstacles represented by circular boundaries, which may not fully capture irregular obstacle geometries or dynamic environments. Third, the analysis is restricted to two-dimensional motion under curvature-constrained kinematics, without explicitly incorporating acceleration, wind disturbances or other vehicle dynamics. Finally, although the model has been evaluated through several numerical experiments and a geographically based case study, further investigation is required to assess its computational scalability in more complex environments and its performance in practical UAV operations. Addressing these limitations would provide opportunities to extend the proposed framework to more realistic and dynamic interception scenarios.
\noindent 
Future work will extend the model to three-dimensional motion, manoeuvring or uncertain targets, and moving obstacles; develop global optimisation or initialisation strategies and real-time replanning; and provide a quantitative comparison with sampling-based planners.

\section{Declarations}

\noindent \textbf{\large {Ethics approval and consent to participate:}} Not applicable.

\noindent \textbf{\large {Consent for publication:}} Not applicable.

\noindent \textbf{\large {Funding:}} 
M. Akter extends her sincere gratitude to the Research $\&$ Publication Cell, University of Chittagong, Bangladesh, for their financial support under Reference No. 674/2024-25/2nd Call/08.

\noindent \textbf{\large {Availability of data and materials:}} Not applicable.


\noindent \textbf{\large {Acknowledgments:}} We appreciate the ChatGPT (3.5) AI tools, which was used to reword, edit, and polish textual content for grammar, spelling, and overall style. 


\newpage
\begin{appendix}{Appendix}
\section{Algorithm for the Proposed Optimal Path} \label{algorithm_target}

\noindent We now present the algorithm developed for our proposed Model \eqref{model_1}, which determines the optimal interception path for a pursuer targeting a moving object in the presence of a single obstacle. 

\begin{algorithm}\label{algorithm}
\begin{description}
\item[Step  1] 
{ \textbf{(Input)}} \\  Define the starting locations for a pursuer and a target as $P(x_P(t_0), y_P(t_0), \theta_P(t_0))$ and $T(x_T(t_0), y_T(t_0), \theta_T(t_0))$, respectively. Assign turning radius $R_a$, curvature $\ds a=\frac{1}{R_a}$, and speed $\mathbf{|V_P|}$ and $\mathbf{|V_T|}$. Set $\ell_i$, $i=1,\dots, 7$, as variables. Set the center of the obstacle $(x_b,y_b)$, its turning radius $R_b\geq R_a$ and the corresponding curvature bound $\ds a_w=\frac{1}{R_b}$.
\item[Step 2] {\textbf{(Find the distance between the pursuer's starting point and the intercept location.)}}\\
Consider the following: left turn angle $\theta_P(t_1)=\theta_P(t_0) + a \ell_1$; right turn angle $\theta_P(t_2) =\theta_P(t_1) - a \ell_2$. Also on the obstacles turning circle, $\theta_P(t_4)=\theta_P(t_2) + a_w \ell_4$ for left turn, ~~$ \theta_P(t_5) =\theta_P(t_4) - a_w \ell_5 $ for right turn. For each of the four path types in \eqref{pathset}, fix the initial turn ($\ell_2=0$ or $\ell_1=0$) and the touring direction $\sigma$ with \eqref{sigmacond}; discard candidates that enter the obstacle turning circle (Remark \ref{rem:inactive}) and keep the shortest one.

To solve Model \eqref{model_1}, find $\ell:=(\bar{\ell}_1, \ldots,  \bar{\ell}_7)$.

Set $\bar{F}$ := $[ \bar{\ell}_1, \ldots, \bar{\ell}_7]$. Set $\bar{T} := [ \bar{t}_1, \ldots, \bar{t}_7]$, where $\bar{t}_i=\bar{\ell}_i/|\mathbf{V_P}|$, $i=1,\ldots,6.$, $\bar{t}_7=\bar{\ell}_7/|\mathbf{V_T}|$ and $f =\ds \sum_{i=1}^7 \bar{\ell}_i$.
 \item[Step 3] {\textbf{(Simulation)}} \\ Same as \cite[Step 3 in Algorithm 5.1]{Akter2023}.

 \item[Step 4] {\textbf{(Output)}} \\
  Calculate $\cup_{i=1}^7\bar{X}_i$.
\end{description}
    
\end{algorithm}
\color{black}
  
\end{appendix}

\begin{thebibliography}{100}

\bibitem{Yong2006}{Yong, C., $\&$ Barth, E. J. (2006)}. Real-time dynamic path planning for Dubins' nonholonomic robot. In Proceedings of the 45th IEEE Conference on Decision and Control, 2418-2423.

\bibitem{Höffmann2024}{Höffmann, M., Patel, S., $\&$ Büskens, C. (2024)}. Optimal guidance track generation for precision agriculture: A review of coverage path planning techniques. Journal of Field Robotics, 41(3), 823-844.

\bibitem{chou2019}{Chou, J. S., Cheng, M. Y., Hsieh, Y. M., Yang, I. T., $\&$ Hsu, H. T. (2019)}. Optimal path planning in real time for dynamic building fire rescue operations using wireless sensors and visual guidance. Automation in construction, 99, 1-17.



\bibitem{Dubins1957} {Dubins, L. E.} (1957). On curves of minimal length with a constraint on average curvature, and with prescribed initial and terminal positions and tangents. American Journal of mathematics, 79(3), 497-516.



\bibitem{Reeds1990} {Reeds, J., $\&$ Shepp, L.} (1990). Optimal paths for a car that goes both forwards and backwards. Pacific journal of mathematics, 145(2), 367-393.

\bibitem{Ayala2014} {Ayala, J., Kirszenblat, D., $\&$ Rubinstein, J. H.} (2014). A geometric approach to shortest bounded curvature paths. arXiv preprint arXiv:1403.4899.


\bibitem{Sussmann1991} {Sussmann, H. J., $\&$ Tang, G.} (1991). Shortest paths for the Reeds-Shepp car: a worked out example of the use of geometric techniques in nonlinear optimal control. Rutgers Center for Systems and Control Technical Report, 10, 1-71.

\bibitem{Boissonnat1994} {Boissonnat, J. D., Cerezo, A., $\&$ Leblond, J.} (1994). Shortest paths of bounded curvature in the plane. Journal of Intelligent and Robotic Systems, 11, 5-20. 



\bibitem{Jha2020} {Jha, B., Chen, Z., $\&$ Shima, T.} (2020). On shortest Dubins path via a circular boundary. Automatica, 121, 109192.

\bibitem{Manyam2019a} {Manyam, S. G., Casbeer, D., Von Moll, A. L., $\&$ Fuchs, Z.} (2019). Shortest Dubins path to a circle. In AIAA scitech 2019 forum, 919.

\bibitem{Chen2019}{Chen, Z., $\&$ Shima, T.} (2019). Shortest Dubins paths through three points. Automatica, 105, 368-375.


\bibitem{Pontryagin1962} {Pontryagin, L. S., Boltyanskii, V.G, Gamkrelidze, R. V., $\&$ Mishchenko, E. F.} (1962) The mathematical theory of optimal processes (Russian), John Wiley and Sons, Interscience Publishers, New York.

\bibitem{Cohen2017}{Cohen, I., Epstein, C., $\&$ Shima, T.} (2017). On the discretized dubins traveling salesman problem. IISE Transactions, 49(2), 238-254.


\bibitem{Madveev2011} {Matveev, A. S., Teimoori, H., $\&$ Savkin, A. V.} (2011). A method for guidance and control of an autonomous vehicle in problems of border patrolling and obstacle avoidance. Automatica, 47(3), 515-524.

\bibitem{Madveev2012} {Matveev, A. S., Teimoori, H., $\&$ Savkin, A. V.} (2012). Method for tracking of environmental level sets by a unicycle-like vehicle. Automatica, 48(9), 2252-2261.

\bibitem{Gopalan2016} {Gopalan, A., Ratnoo, A., $\&$ Ghose, D.} (2016). Time-optimal guidance for lateral interception of moving targets. Journal of Guidance, Control, and Dynamics, 39(3), 510-525. 

\bibitem{Gopalan2017} {Gopalan, A., Ratnoo, A., $\&$ Ghose, D.} (2017). Generalized time-optimal impact-angle-constrained interception of moving targets. Journal of Guidance, Control, and Dynamics, 40(8), 2115-2120.

\bibitem{Manyam2019} {Manyam, S. G., Casbeer, D., Von Moll, A., \& Fuchs, Z.} (2019). Optimal dubins paths to intercept a moving target on a circle. In 2019 American Control Conference (ACC), 828-834.

\bibitem{Buzikov2021} {Buzikov, M. E., $\&$ Galyaev, A. A.} (2021). Time-optimal interception of a moving target by a Dubins car. Automation and Remote Control, 82, 745-758.

\bibitem{Looker2008} {Looker, J. R.} (2008). Minimum Paths to Interception of a Moving Target when Constrained by Turning Radius.


\bibitem{Meyer2015} {Meyer, Y., Isaiah, P., $\&$ Shima, T.} (2015). On Dubins paths to intercept a moving target. Automatica, 53, 256-263.

\bibitem{Triharminto2014} {Triharminto, H. H., Prabuwono, A. S., Sulaiman, R., Setiawan, N. A., $\&$ Adji, T. B.} (2014). Turning strategy based on modified Dubins algorithm for unmanned aerial vehicle dynamic path planning. In Proceedings of the 2013 International Conference on Advances in Intelligent Systems in Bioinformatics, Atlantis Press, 82-86.

\bibitem{Yang2002}{Yang, G., $\&$ Kapila, V.} (2002, December). Optimal path planning for unmanned air vehicles with kinematic and tactical constraints. In Proceedings of the 41st IEEE Conference on Decision and Control, 2002. Vol. 2, 1301-1306.

\bibitem{Kaya2017} {Kaya, C. Y.} (2017). Markov–Dubins path via optimal control theory. Computational Optimization and Applications, 68(3), 719-747.

\bibitem{Kaya2019} {Kaya, C. Y.} (2019). Markov–Dubins interpolating curves. Computational Optimization and Applications, 73(2), 647-677.

\bibitem{Zheng2021} {Zheng, Y., Chen, Z., Shao, X., $\&$ Zhao, W.} (2021). Time-optimal guidance for intercepting moving targets by dubins vehicles. Automatica, 128, 109557.

\bibitem{Zheng2022}{Zheng, Y., Zheng, C. H. E. N., Xueming, S. H. A. O., \& Wenjie, Z. H. A. O.} (2022). Time-optimal guidance for intercepting moving targets with impact-angle constraints. Chinese Journal of Aeronautics, 35(7), 157-167.

\bibitem{Forkan2022}{Forkan, M., Rizvi, M. M., $\&$ Chowdhury, M. A. M.} (2022). Optimal path planning of Unmanned Aerial Vehicles (UAVs) for targets touring: Geometric and arc parameterization approaches. Plos one, 17(10), e0276105.

\bibitem{Ferrari2009}{Ferrari, S., Fierro, R., Perteet, B., Cai, C., $\&$ Baumgartner, K.} (2009). A geometric optimization approach to detecting and intercepting dynamic targets using a mobile sensor network. SIAM journal on control and optimization, 48(1), 292-320.

\bibitem{Hou2023}{Hou, L., Luo, H., Shi, H., Shin, H. S., $\&$ He, S.} (2023). An Optimal Geometrical Guidance Law for Impact Time and Angle Control. IEEE Transactions on Aerospace and Electronic Systems.



\bibitem{Cao2015}{Cao, Y.} (2015). UAV circumnavigating an unknown target under a GPS-denied environment with range-only measurements. Automatica, 55, 150-158.







\bibitem{Parlangeli2024}{Parlangeli, G., De Palma, D., $\&$ Attanasi, R.} (2024). A novel approach for 3PDP and real-time via point path planning of Dubins’ vehicles in marine applications. Control Engineering Practice, 144, 105814.










































\bibitem{Fourer2003} {Fourer, R., Gay, D. M., $\&$ Kernighan, B. W.} (2003). AMPL: A Modeling Language for Mathematical Programming, Second Edition, Brooks/Cole Publishing Company / Cengage Learning.

\bibitem{Wachter2006} {Wächter, A., $\&$ Biegler, L. T.} (2006). On the implementation of an interior-point filter line-search algorithm for large-scale nonlinear programming. Mathematical programming, 106, 25-57.




\bibitem{Akter2023}{Akter, M., Rizvi, M. M., $\&$ Forkan, M. (2023)}. A New Algorithm for Optimizing Dubins Paths to Intercept a Moving Target. arXiv preprint arXiv:2310.17856.




\end{thebibliography}
\end{document}